\documentclass[12pt]{article}

\usepackage{latexsym,amssymb}

\def\tr{^{\rm T}}
\def\Real{\mathbb R}

\def\dst{\displaystyle}

\def\qmx#1{\left(\matrix{#1}\right)}

\def\beeq#1{\begin{equation}{#1}\end{equation}}

\def\ba{\begin{array}}
\def\ea{\end{array}}
\def\eqa{\begin{eqnarray}}
\def\eqe{\end{eqnarray}}

\newtheorem{proposition}{Proposition}

\newtheorem{lemma}{Lemma}

\newenvironment{proof}{\medskip\noindent{\it Proof. }}{ \medskip}
\newenvironment{remark}{\medskip\noindent{\it Remark. }}{
\medskip}

\begin{document}

\title{A New Approach to Adaptive Nonlinear Regulation \thanks{This work was partially
supported by NSF under grant ECS-0314004, by ONR under grant
N00014-03-1-0314. Corresponding Author: Dr. Lorenzo Marconi,
email: lmarconi@deis.unibo.it, tel. 0039 051 2093788, fax. 0039
051 2093073}}

\author{F. Delli Priscoli $^{\dag}$,  L. Marconi $^{\circ}$,
A.Isidori $^{\dag \circ\ddag
}$}
\date{\today}

\maketitle
\begin{center}

\small $^{\dag}$Dipartimento di Informatica e Sistemistica,
 Universit\`{a} di Roma ``La Sapienza'', \\00184 Rome, ITALY.

$^{\circ}$ C.A.SY. -- Dipartimento di Elettronica, Informatica e
Sistemistica, University of Bologna, \\40136 Bologna, ITALY.

$^{\ddag}$Department of Electrical and Systems Engineering,
Washington University, \\St. Louis, MO 63130.
\smallskip
\normalsize
\end{center}
\maketitle \vspace{1cm}
\begin{center}
{\em \large Paper submitted to SIAM Journal On Control and
Optimization}
\end{center}
\vskip0.5in
\begin{abstract}
This paper shows how the theory of adaptive observers can be
effectively used in the design  internal models for nonlinear
output regulation. The main result obtained in this way is a new
method for the synthesis of adaptive internal models which
substantially enhances the existing theory of adaptive output
regulation, by allowing nonlinear internal models and more general
classes of controlled plants.
\end{abstract}


\section{Introduction}\label{sec1}

The problem of controlling the output of a system so as to achieve
asymptotic tracking of prescribed trajectories and/or asymptotic
rejection of disturbances is a central problem in control theory.
There are essentially three different possibilities to approach
the problem: tracking by dynamic inversion, adaptive tracking,
tracking via internal models. Tracking by dynamic inversion
consists in computing a precise initial state and a precise
control input (or equivalently a reference trajectory of the
state), such that, if the system is accordingly initialized and
driven, its output exactly reproduces the reference signal. The
computation of such control input, though, requires  ``perfect
knowledge" of the entire trajectory to be tracked as well as
``perfect knowledge" the model of the controlled plant. Thus, this
type of approach is not suited in the presence of large
uncertainties on plant parameters as well as on the reference
signal. Adaptive tracking can successfully handle parameter
uncertainties, but it still presupposes the knowledge of the
entire trajectory which is to be tracked (to be used in the design
of the adaptation algorithm) and therefore this approach is not
suited in the problem of tracking unknown trajectories.
Internal-model-based tracking on the other hand, is able to handle
simultaneously uncertainties in plant parameters as well as in the
trajectory which is to be tracked. It has been proven that, if the
trajectory to be tracked belongs to the set of all trajectories
generated by some fixed dynamical system, a controller which
incorporates an internal model of such a system is able to secure
asymptotic decay to zero of the tracking error for every possible
trajectory in this set and does it robustly with respect to
parameter uncertainties. This is in sharp contrast with the two
approaches mentioned above, where in lieu of the assumption that a
signal is within a class of signals generated by an exogenous
system, one instead needs to assume complete knowledge of the
past, present and future time history of the trajectory to be
tracked. It is for this reason that the internal-model-based
approach seems to be the best suited in problems of tracking of
unknown reference trajectories or rejecting unknown disturbances.

A generalized problem of tracking and asymptotic disturbance
rejection is usually cast as follows. A nonlinear system is given,
modelled by equations of the form
\[
\begin{array}{rcl}\dot x &=& f(x,u,v)\\y&=&k(x,v)\\ e
&=&h(x,v)
\end{array}
\]
with state $x$, control input $u$, measured output $y$, regulated
output $e$. In this system, $v$ is an exogenous input, which
represents actual {\em disturbances} as well as {\em commands} to
be followed, and it is assumed that, as a function time, $v(t)$
can be seen as generated by a separate autonomous dynamical
system, called the {\em exosystem}. Generally speaking, the
problem of tracking and asymptotic disturbance rejection
(sometimes also referred to as the generalized {\em servomechanism
problem} or the {\em output regulation problem}) is to design a
controller so as to obtain a closed-loop system in which:
\begin{itemize} \item all trajectories are bounded, and \item the
regulated output $e(t)$ asymptotically decays to $0$ as $t\to
\infty$ .\end{itemize} The peculiar aspect of this design problem
is the characterization of the class of all possible exogenous
inputs (disturbances as well commands) as the set of all possible
solutions of a fixed (finite-dimensional) differential equation.
This can be seen as an intermediate choice sitting between two
extremes: the (pessimistic) case in which the design is required
to obtain certain goals in the presence of the {\em worst
possible} exogenous input and the (optimistic) case in which the
controller is assumed {\em to have access} to the exogenous input
$v$. In this design problem, the controller does not have access
to the exogenous input in real time, but the latter is restricted
to range  over a ``finite dimensional" set of functions (such as
the set of solutions of a fixed differential equation). The vector
$v$ may include constant uncertain parameters, which have a
trivial dynamics and hence can be viewed as solutions of a
(trivial) differential equation. In other words, in this setting,
{\em any source of uncertainty} (about an actual disturbance
affecting the system, about an actual trajectory to be tracked,
about any unknown constant parameter in the plant or about any
unknown constant parameter in the exosystem itself) is treated as
{\em uncertainty in the initial condition} of a fixed autonomous
finite dimensional dynamical system, which is then seen as source
of all possible, constant as well time-varying, uncertainties.

For linear multivariable systems this problem was addressed in
very elegant geometric terms by Davison, Francis, Wonham \cite{Da,
FW76, Fr} and others.  A nonlinear enhancement of this theory,
which uses a combination of geometry and nonlinear dynamical
systems theory, was presented in \cite{IB, HR, HL, BDIK} in the
context of solving the problem near an equilibrium, in the
presence of exogenous signals which were produced by a Poisson
stable system. In particular, Huang showed how, by appropriately
designing the internal model, the controlled output could be
steered to zero in spite of plant parameter uncertainties, thus
extending to the nonlinear setting one of the most remarkable
features of internal-model-based design for linear systems. Under
suitable hypotheses, the (local) design methods presented in these
works have been extended \cite{Kh2, IA, SIM99a} to the case of
arbitrary large (but compact) sets of initial data. A substantial
limitation of classical internal model-based control (for linear
as well as nonlinear systems) is the sensitivity to parameters
uncertainties in the exosystem. This limitation, though, was later
addressed and solved under convenient hypotheses in the paper
\cite{SIM}, where the possibility of using techniques of adaptive
control to cope with unknown parameters in the exosystem was
successfully demonstrated.

In the recent paper \cite{BI03}, the problem in question has been
posed in more general terms, not tied, as all previous
contributions were, to the existence of a privileged equilibrium
point about which the (local as well semi-global) analysis was
conducted. The more general foundations laid in this way make it
possible to overcome certain restrictions of the earlier theory,
notably the assumption that the controlled plant has an
asymptotically stable zero-dynamics, which is replaced by the
substantially weaker hypothesis that the latter possess a compact
attractor. Another major enhancement of this newer approach is a
systematic method for the design of nonlinear internal models (see
\cite{BI03bis}). The presence of parametric uncertainties in the
exosystem, however, is not explicitly addressed in these works.

The purpose of the present paper is to show how the problem of
handling {\em parametric uncertainties in the exosystem} can be
successfully addressed by means of a new approach which reposes,
on one hand, on the general non-equilibrium theory developed in
\cite{BI03} and, on the other hand, on the theory of adaptive
observers for nonlinear system pioneered in \cite{BaGe} and
\cite{MaTo92}. The result obtained in this way is a totally new
method for the synthesis of adaptive internal models which
substantially extends the adaptive regulation theory presented in
\cite{SIM}, by allowing nonlinear internal models and more general
classes of controlled plants.

\section{Output regulation and limit sets}\label{basicb}

The purpose of output regulation is to obtain a closed-loop system
in which all trajectories with initial conditions in a fixed (but
otherwise arbitrary) compact set are bounded and the regulated
output converges to zero as time tends to infinity. As shown in
\cite{BI03}, intimately associated with this problem is the notion
of {\em limit set}  of a given bounded set of initial conditions.
For convenience of the reader, the notion in question is
summarized as follows.

Consider an {\em autonomous} ordinary differential equation
\beeq{\label{uno}\dot x = f(x)}  in which $x\in \Real^n$, $t\in
\Real$, with $f(x)$ a locally Lipschitz function. Let
\[\ba{rccl}
\phi  :&(t,x)&\mapsto&\phi(t,x) \ea\] define the flow of
(\ref{uno}). Suppose the flow is forward complete.  The
$\omega$-{\em limit set} of a subset $B \subset \Real^n$, written
$\omega(B)$, is the totality of all points $x\in\Real^n$ for which
there exists a sequence of pairs $(x_k,t_k)$, with $x_k\in B$ and
$t_k\to \infty$ as $k\to \infty$, such that
\[
\lim_{k\to \infty}\phi(t_k,x_k)=x\,. \] In case $B=\{x_0\}$ the
set thus defined, $\omega(x_0)$, is precisely  the $\omega$-limit
set, as defined by Birkhoff, of the point $x_0$. Note that, in
general
\[
\bigcup_{x_0\in B}\omega(x_0) \subset \omega(B)\,, \] but the
equality may not hold.

It is well-known that $\phi(t,x_0)$, if is bounded in positive
time, the set $\omega(x_0)$ is non-empty, compact, invariant, and
\[
\lim_{t \to \infty}{\rm dist}(\phi(t,x_0),\omega(x_0))=0\,.
\]
If $B$ is not just the singleton $\{x_0\}$, the following more
general property holds. Recall that a set $A$ is said to {\em
uniformly attract} a set $B$ under the flow of (\ref{uno}) if for
every $\varepsilon
>0$ there exists a time $\bar t$ such that
\[
{\rm dist}(\phi(t,x),A)\le \varepsilon, \qquad \mbox{for all $t\ge
\bar t$ and for all $x\in B$.}
\]
Then the following holds (see \cite[page 8]{Hale91}).

\begin{lemma}\label{LM0} If $B$ is a nonempty bounded set for which there is
a compact set $J$ which uniformly attracts $B$ (thus, in
particular, if $B$ is any nonempty bounded set whose positive
orbit has a bounded closure), then $\omega(B)$ is nonempty,
compact, invariant and uniformly attracts $B$. Moreover, if
$\omega(B)\in {\rm int}(B)$, then $\omega(B)$ is stable in the
sense of Lyapunov.
\end{lemma}

\section{Class of systems and main assumptions}\label{SecPA}
In this paper we discuss the design of output regulators for
nonlinear systems modelled by equations of the form
 \beeq{\label{pla1}\ba{rcl} \dot z &=& f_0(\varrho,w,z) +
 f_1(\varrho,w, z,e_1)e_1\\
 \dot e_1 &=& e_2\\
 &\vdots&\\
 \dot e_{r-1} &=& e_r\\
 \dot e_r &=& q(\varrho,w,z,e_1, \ldots,e_r) + u\\
 e &=& e_1\\
 y &=& \mbox{col}(e_1,\ldots,e_r)\,,
 \ea}
 with state $(z,e_1,\ldots,e_r)\in \Real^n \times \Real^r$,
control input $u\in \Real$, regulated output $e\in \Real$,
measured output $y\in \Real^r$, in which the exogenous
(disturbance) input $w\in \Real^s$ is generated by an exosystem
\beeq{\label{exo} \dot w = s(\varrho,w)\,.} In this model,
$\varrho \in \Real^p$ is a vector of {\em constant} uncertain
parameters, ranging over a fixed compact set $P$. The vector
$\varrho$ is the aggregate of a finite set of uncertain parameters
affecting the controlled plant and another, possibly different,
set of uncertain parameters affecting the exosystem. These
parameters may be regarded as ``\,trivial components" of an
``\,augmented" exogenous input, but for the sake of clarity, and
also consistency with some of the earlier literature, their role
will be kept separate. Occasionally, throughout the paper, the
``\,augmented" exosystem
 \beeq{\label{exo1}
 \ba{rcl}
 \dot \varrho&=&0\\
 \dot w &=& s(\varrho,w)
 \ea
 }
will be rewritten in more compact form as \beeq{\label{exobold}
 \dot {\bf w} = {\bf s}({\bf w})\,,} where ${\bf w} = {\rm
 col}(\varrho, w)$.

The functions $f_0(\cdot), f_1(\cdot), q(\cdot), s(\cdot)$ in
(\ref{pla1}) and (\ref{exo1}) are assumed to be at least
continuously differentiable. The initial conditions of
(\ref{pla1}) range on a set $Z\times E$, in which $Z$ is a fixed
{\em compact} subset of $\Real^n$ and $E=\{(e_1,\ldots,e_r)\in
\Real^r:|e_i|\le c\}$, with $c$ a fixed number. The initial
conditions of the exosystem (\ref{exobold}) range on a {\em
compact} subset ${\bf W}$ of $\Real^p \times \Real^s $. In this
framework the problem of output regulation is to design an output
feedback regulator of the form
 \[
 \ba{rcl}
 \dot \zeta &=& \varphi(\zeta, y)\\
 u &=& \gamma(\zeta,y)
 \ea
 \]
 such that {\em for all initial conditions ${\bf w}(0) \in
 {\bf W}$ and $(z(0),e_1(0), \ldots, e_r(0)) \in Z \times E$ the
 trajectories of the closed-loop system are bounded and $\lim_{t \rightarrow \infty}
 e(t)=0$.}\\

Augmenting ({\ref{pla1}) with (\ref{exo1}) yields a system which,
viewing $u$ as input and $e$ as output, has relative degree $r$.
The associated ``augmented" zero dynamics, which is forced by the
control
 \beeq{ \label{friend1}  c(\varrho,
 w, z) = -q(\varrho,w,z,0,\ldots, 0)\,,
 }
is given by
 \beeq{\label{zerodynpla}\ba{rcl}
 \dot \varrho &=&0\\[2mm]
 \dot w &=& s(\varrho,w)\\[2mm]
 \dot z &=& f_0(\varrho,w,z) \,.\ea
 }

Occasionally, throughout the paper, we will find it convenient to
rewrite the latter in more compact form as
 \beeq{\label{zerobold}
 \dot {\bf z} = {\bf f}_0({\bf z})\,,}
 having set
 ${\bf z}={\rm col}(\varrho,w,z)$.
 Accordingly, we set ${\bf Z}= {\bf W}\times Z$ and, with a mild abuse of notation, we
 replace $c(\varrho,
 w, z)$ by $c({\bf z})$ in (\ref{friend1}).

In what follows, we retain three of the basic assumptions that
were introduced in \cite{BI03} and express certain properties of
the augmented zero dynamics (\ref{zerodynpla}). The assumptions in
question are the following ones:

\medskip\noindent {\em Assumption (i)}\,: the set ${\bf W}$ is a differential submanifold (with
boundary) of $\Real^{p} \times \Real^s$, and ${\bf W}$ is
invariant for (\ref{exobold}). $\triangleleft$

\medskip\noindent {\em Assumption (ii)}\,:
there exists a compact subset $\cal Z$ of ${\bf W} \times \Real^n$
which contains the positive orbit of the set ${\bf Z}$ under the
flow of (\ref{zerobold}), and $\omega({\bf Z})$ is a differential
submanifold (with boundary) of ${\bf W} \times \Real^n$. Moreover
there exists a number $d_1>0$ such that
 \[
 {\bf z} \in
{\bf W} \times \Real^n\,, \quad \mbox{dist}({\bf z}, \omega({\bf
Z})) \leq d_1 \qquad \Rightarrow \qquad {\bf z} \in {\bf Z}\,.
\quad \triangleleft
\]

\medskip As a remark on the above hypotheses, note that, since the
positive orbit of the set $\bf Z$ under the flow of
(\ref{zerobold}) is bounded, the set $\omega({\bf Z})$, namely the
$\omega$-limit set of ${\bf Z}$ under the flow of
(\ref{zerobold}), is a nonempty, compact and invariant subset of
${\bf W} \times \Real^n$ which uniformly attracts all trajectories
of (\ref{zerobold}) with initial conditions in ${\bf Z}$. It can
also be shown (as in \cite{BI03}) that  for every ${\bf w} \in \bf
W$ there is $z \in \Real^n$ such that $({\bf w}, z) \in
\omega({\bf Z})$. In what follows, for convenience, the set
$\omega({\bf Z})$ will be simply denoted as ${\cal A}_0$.

The last condition in assumption (ii) implies that ${\cal A}_0$ is
stable in the sense of Lyapunov. The next hypothesis, which will
be used in the last part of the paper, is that the set ${\cal
A}_0$ is locally exponentially attractive.

\medskip\noindent {\em Assumption (iii)}\,: There exist
$M \geq 1$, $a>0$ and $d_2 \leq d_1$ such that
\[
{\bf z}_0 \in {\bf W} \times \Real^n\,, \quad \mbox{dist}({\bf
z}_0, {\cal A}_0) \leq d_2 \qquad \Rightarrow \qquad
 \mbox{dist}({\bf z}(t,{\bf z}_0), \, {\cal A}_0) \leq M e^{- a t}
 \mbox{dist}({\bf z}_0 , \, {\cal A}_0)
\]
in which ${\bf z}(t,{\bf z}_0)$ denotes the solution of
(\ref{zerobold}) passing through ${\bf z}_0$ at time $t=0$.
$\triangleleft$

\medskip
The results presented in \cite{BI03}, as essentially all previous
results on output regulation, relied upon the hypothesis that the
set of all ``feed-forward inputs capable to secure perfect
tracking" (that, is, the set of inputs of the form $u(t) = c({\bf
z}(t))$, with ${\bf z}(t)$ a trajectory of the restriction of
(\ref{zerobold}) to ${\cal A}_0$) could be seen as a subset of the
set of outputs of a suitable linear system. The system in question
was used to construct a (linear, as a matter of fact) internal
model. This assumption was weakened in \cite{BI03bis}, where a
general method for the construction of fully nonlinear internal
models was presented, but the method in question did not allow for
the presence of uncertain parameters in the exosystem. In this
paper we introduce a different kind of hypothesis, leading to a
somewhat more restricted class of internal models, but which -- in
return -- allows for uncertain parameters in the exosystem.

\medskip\noindent {\em Assumption (iv)}\,:
 there exist a positive integer $d$, a $C^1$ map
 \[
 \ba{lrcl}
 \tau \quad : \quad &  {\cal Z} &\rightarrow& \Real^d\\
              & {\bf z} &\mapsto& \tau({\bf z})\,,
 \ea
 \]
 a $C^0$ map
 \[
 \ba{lrcl}
 \theta \quad : \quad &  {P} &\rightarrow& \Real^q\\
              & \varrho &\mapsto& \theta(\varrho)\,,
 \ea
 \]
 an observable pair $(A,C) \in \Real^{d \times d} \times
 \Real^{1 \times d}$, and two $C^1$ maps $\phi: \Real \to \Real^d$ and
 $\Omega: \Real \to \Real^{d \times q}$ such that the following
 identities (which we call
 {\em immersion property})
 \beeq{ \label{adaptimmers1}
 {\partial \tau \over \partial {\bf z}} \, {\bf f}_0({\bf z})
   =A\, \tau({\bf z}) + \phi(C\tau({\bf z}))
  + \Omega(C\tau({\bf z}))\, \theta(\varrho)
 }
 \beeq{\label{adaptimmers0}
 c({\bf z}) = C\, \tau({\bf z})\,
 }
 hold for all ${\bf z}\in {\cal A}_0$, $\varrho \in P$. $\triangleleft$

 \begin{remark}
 Without loss of generality (see \cite[page 208]{MaToBook}),
 we can assume throughout that the matrices $A$ and $C$ in
 (\ref{adaptim}) have the form
 \[
 A =\qmx{0 & 1 & 0 &\cdots &0\cr 0 & 0 & 1 &\cdots &0\cr \cdot&
 \cdot & \cdot & \cdots &\cdot \cr 0 & 0 & 0 &\cdots &1\cr 0 & 0 &
 0& \cdots& 0\cr},\qquad C = \qmx{1 & 0 & 0 &\cdots &0\cr}.
 \]
 Furthermore, note that since the maps
 $\Omega(\cdot)$ and $\phi(\cdot)$ are continuously differentiable and the relations
 (\ref{adaptimmers1}) -- (\ref{adaptimmers0})
 are supposed to hold over the compact set ${\cal A}_0$, it can be assumed without loss of
 generality that functions $\phi(\cdot)$ and $\Omega(\cdot)$ have compact support.
 This being the case, the functions in question can be assumed
 {\em globally Lipschitz}, i.e. there exist $L_\phi$ and $L_\Omega$ such that
 \[
 |\phi(s_1)-\phi(s_2)| \le L_\phi|s_1-s_2|, \qquad
 |\Omega(s_1)-\Omega(s_2)| \le L_\Omega|s_1-s_2|,
 \]
 for all $s_1, s_2$. $\triangleleft$
 \end{remark}

 \begin{remark} Note that Assumption (iv) can be rephrased by saying
 that for each initial condition ${\bf z}(0) \in {\cal A}_0$ of (\ref{zerobold}),
 there is a pair
 $\xi(0),\theta$ such that the control input
 $u(t)=c({\bf z(t)})$ (which is the unique input capable of
 keeping $e(t)$ identically at zero) can be seen as output of a system
 of the form
 \beeq{ \label{adaptim}
 \ba{rcl}\dot \xi &=& A\xi + \phi(y) +
 \Omega(y)\theta\\[2mm]\dot \theta &=& 0\\[2mm]
 y&=&C\xi\,. \qquad\qquad\triangleleft\ea
 }
 \end{remark}

 In the remaining part of this section we show that {\em there is no loss of
 generality} in addressing the simpler case in which the
 relative degree of (\ref{pla1}) is $r=1$.
 As a matter of fact consider the change of variable
 \[
 e_r \; \mapsto\;
 \tilde e := e_r + g^{r-1} a_0 e_1 + g^{r-2} a_1 e_2 + \ldots + g
 a_{r-2} e_{r-1}
 \]
 where $g$ is a {\em positive} design parameter and
 $a_i$, $i=0,\ldots, r-2$\,, are such that
 all roots of the polynomial  $\lambda^{r-1}+a_{r-2} \lambda^{r-1} + \ldots + a_{1} \lambda + a_0
 =0$ have negative real part. This changes system (\ref{pla1})
 into a system of the form
\beeq{\label{plaaux}\ba{rcl} \dot {\tilde z} &=& \tilde
f_0(\varrho,w,\tilde z) +
 \tilde f_1(\varrho,w, \tilde z,\tilde e)\tilde e\\
\dot {\tilde e} &=& \tilde q(\varrho,w,\tilde z, \tilde e,g) + u
 \ea}
in which
\[
\tilde z = {\rm col}(z,e_1,\ldots, e_{r-1})
\]
\[
\tilde f_0(\varrho,w,\tilde z) =
\qmx{f_0(\varrho,w,z)+f_1(\varrho,w,z,e_1)e_1\cr e_2 \cr \cdots
\cr e_{r-1}\cr- g^{r-1} a_0 e_1 - g^{r-2} a_1 e_2 - \ldots - g
 a_{r-2} e_{r-1}\cr}\,\qquad \tilde f_1(\varrho,w,\tilde z,\tilde e) =\qmx{0\cr
 0\cr \cdots \cr 0\cr 1\cr}\] and
 \beeq{\label{tildeqq}
 \ba{rcl}\tilde q(\varrho,w,\tilde z,\tilde e,g) &=& q(\varrho, w, z, e_1, \ldots
 ,e_r)- g^{r-1}a_0e_2 - \cdots -g^2a_{r-3}e_{r-1}\\[2mm]&&\qquad-\; ga_{r-2}[\tilde e- g^{r-1} a_0 e_1 - g^{r-2} a_1 e_2 - \ldots - g
 a_{r-2} e_{r-1}]\,.\ea
 }
 Let the initial conditions of (\ref{plaaux}}) range on a set
 of
 the form $Z\times Z_e \times \tilde E$, in which $Z_e = \{(e_1,\ldots,e_{r-1}:|e_i|\le
 c\}$ and $\tilde E=\{\tilde e: |\tilde e| \le \tilde c\}$ with
 \[
 \tilde c \ge  (1 + g^{r-1} a_0  + g^{r-2} a_1  + \ldots + g
 a_{r-2} )c\] (note the
dependence on the
 choice of the $a_i$'s and of $g$).

 Let system (\ref{plaaux}) be augmented with (\ref{exo1}) and consider a
 regulation problem with {\em regulated} output $\tilde e$ and {\em
 measured}
 output $\tilde y = \tilde e$. The system, viewed as a system with input $u$ and output
$\tilde e$, has relative
  degree $1$ and its zero
  dynamics, forced by the control
 \beeq{ \label{friend2}
 \tilde c (\varrho, w, \tilde z) = -\tilde q(\varrho,w,\tilde z,
 0,g)\,,
 }
 is given by
 \beeq{\label{zdaux}\ba{rcl}
 \dot \varrho&=&0\\
 \dot w &=& s(\varrho,w)\\
 \dot {\tilde z} &=& \tilde f_0(\varrho,w,\tilde z) \,.
 \ea}
 Consistently with the notation used for (\ref{zerodynpla}), the latter can be rewritten in more succinct
 form as
 \beeq{\label{zeroboldtilde}
 \dot {\tilde {\bf z}} = \tilde {\bf f}_0(\tilde {\bf z}) \qquad
 \mbox{with} \qquad
 \tilde {\bf z} = \mbox{col}(\varrho, w, \tilde z )
 \,.
 }

 Suppose that a controller of the form
 \beeq{\label{dummycont}
 \ba{rcl}
 \dot {\zeta} &=& \varphi(\zeta, \tilde y)\\
 u &=& \gamma(\zeta, \tilde y)
 \ea
 }
 has been found which solves the problem of output regulation
 thus defined.  Then, it is immediate to realize that the
 controller
\beeq{ \label{regrd}
 \ba{rcl}
 \dot {\zeta} &=& \varphi(\zeta, e_r + g^{\, r-1}\, a_0\, e_1 + g^{\, r-2} \,a_1 \,
 e_2 + \ldots +  g \, a_{r-2} \,e_{r-1})\\
 u &=& \gamma(\zeta, e_r +  g^{\, r-1} \,a_0 \, e_1 +  g^{\, r-2} \, a_1 \, e_2 + \ldots +
  g \, a_{r-2} \, e_{r-1})
 \ea
 }
 solves the problem of output regulation for the original plant (\ref{pla1}).
 To this end note, first of all, that  (\ref{regrd}) is an admissible controller for (\ref{pla1}),
 because it is driven only by the components $e_1, \ldots, e_r$ of the measured output $y$ of
 (\ref{pla1}). Trivially, the composition of (\ref{pla1}) with
 (\ref{regrd}) differs from the composition of (\ref{plaaux}) with
 (\ref{dummycont}) only by a linear change of coordinates, and for any initial state of (\ref{pla1}) in $Z\times
 E$, the corresponding initial state of (\ref{plaaux}) is in $Z\times Z_e\times
 \tilde E$. Thus all trajectories of (\ref{pla1}), controlled by (\ref{regrd}), with initial conditions in
 $Z\times E$ are bounded. The trajectories in question are such
 that $\lim_{t \to \infty} \tilde e(t) =0$. But since
 \[\ba{rcl}
\dot e_1 &=& e_2\\
 &\vdots&\\
 \dot e_{r-1} &=& - (g^{r-1} a_0 e_1 + g^{r-2} a_1 e_2 + \ldots + g
 a_{r-2} e_{r-1})+\tilde e\ea\]
 and the $a_i$'s are coefficients of a Hurwitz polynomial, it is
 readily concluded that also $\lim_{t \to \infty} e_1(t)
 =0$. Therefore (\ref{regrd}) solves the problem of output regulation for the system
(\ref{pla1})
 .\\

 In the light of these considerations, what is left to show in
 order to prove the desired claim (namely the fact that there is no loss of generality in
 addressing the problem for systems having relative degree 1) is that
 the zero dynamics (\ref{zdaux}) and the associated map
 (\ref{friend2}) inherit, from  (\ref{zerodynpla})
 and (\ref{friend1}), the appropriate properties which make the solution of the problem of
 output regulation possible.
 Specifically, we will prove that if
  (\ref{zerodynpla}) and (\ref{friend1}) satisfy assumptions (i)-(iv)
 above, then (\ref{zdaux}) and (\ref{friend2})
 satisfy an {\em identical} set of assumptions, {\em provided that} the parameter $g$ is chosen
 sufficiently large.
 This is formalized in the next Lemma.

 \begin{lemma}\label{lemmard}
 Suppose that assumptions (i)-(iv) hold for
 (\ref{zerodynpla}) and (\ref{friend1}). Set $\tilde {\bf Z} = {\bf W}\times Z\times Z_e$.
 Then there exists $g^\star>0$ such that
 for all fixed $g \geq g^\star$ the following hold:

 \medskip\noindent {${(ii)'}$}
there exists a compact subset $\tilde {\cal Z}$ of ${\bf W} \times
\Real^n \times \Real^{r-1}$ which contains the positive orbit of
the set $\tilde {\bf Z}$ under the flow of (\ref{zeroboldtilde}),
and $\tilde {\cal A}_0:=\omega(\tilde {\bf Z})$ is a differential
submanifold (with boundary) of ${\bf W} \times \Real^n \times
\Real^{r-1}$. Moreover there exists a number $\tilde d_1>0$ such
that
\[
 \tilde {\bf z} \in
{\bf W} \times \Real^n \times \Real^{r-1}\,, \quad \mbox{\em
dist}(\tilde {\bf z}, \tilde {\cal A}_0) \leq \tilde d_1 \qquad
\Rightarrow \qquad \tilde {\bf z} \in \tilde{\bf Z}\,.
\]

 \medskip \noindent {${(iii)'}$} there exist
$\tilde M \geq 1$, $\tilde a>0$ and $\tilde d_2 \leq \tilde d_1$
such that
\[
 \tilde {\bf z}_0 \in {\bf W} \times \Real^n \times \Real^{r-1}\,,
 \quad \mbox{\em dist}(\tilde {\bf
 z}_0, \tilde {\cal A}_0) \leq \tilde d_2 \qquad \Rightarrow \qquad
 \mbox{\em dist}(\tilde {\bf z}(t,\tilde {\bf z}_0), \, \tilde {\cal A}_0) \leq \tilde
 M e^{- \tilde a t}
 \mbox{\em dist}(\tilde {\bf z}_0 , \, \tilde {\cal A}_0)
\]
in which $\tilde {\bf z}(t,\tilde {\bf z}_0)$ denotes the solution
of (\ref{zeroboldtilde}) passing through $\tilde {\bf z}_0$ at
time $t=0$.

 \medskip \noindent {${(iv)'}$} there exist a $C^1$ map
  \[
 \ba{lccl}
 \tilde \tau \quad : \quad &  \tilde {\cal Z} &\rightarrow& \Real^d\\
              & \tilde {\bf z}
              &\mapsto& \tilde \tau(\tilde{\bf z})
 \ea
 \]
 such that the immersion property\,\footnote{As above, with a mild abuse of notation we rewrite $\tilde
 c(\varrho, w, z, e_1, \ldots, e_{r-1})$ as $\tilde c(\tilde {\bf
 z})$.}

  \[
  {\displaystyle \partial {\tilde \tau} \over
  \displaystyle \partial \tilde {\bf z} } \,\tilde {\bf f}_0(\tilde {\bf z})
  = A \tilde \tau(\tilde {\bf z}) + \phi(C \tilde\tau(\tilde{\bf z}))
 + \Omega(C\tilde \tau(\tilde{\bf z}))\, \theta(\varrho)
 \]
 \[
 \tilde c (\tilde {\bf z}) = C \tilde \tau(\tilde {\bf z})
 \]
 holds for all $\tilde {\bf z} \in \tilde {\cal A}_0$, and $\varrho \in P$.
 \end{lemma}

\begin{proof}
Consider the change of variable
 \[
  x = \qmx{x_1 \cr x_2 \cr \vdots \cr x_{r-1}\cr}: = D_g^{-1} \qmx{e_1 \cr e_2 \cr \vdots
  \cr e_{r-1}\cr} \qquad
  \mbox{with} \qquad D_g=\left (\ba{cccc}
  1 & 0 & \ldots & 0\\
  0 & {\dst g} & \ldots & 0\\
  \vdots & \vdots & \ddots& \vdots\\
  0 & 0 & \ldots & {\dst g^{r-2}}
  \ea   \right )
 \]
 which transforms system (\ref{zdaux}) into
 \beeq{\label{sysappb}
 \ba{rcl}
 \dot z &=& f_0(\varrho, w,z) + f_1(\varrho, w,z, x_1)
 x_1\\
 \dot x &=& g A x
 \ea
 }
 where $A$ is a Hurwitz matrix.
 Note that if $g>1$, which we can assume without loss of generality,
 $(e_1(0), \ldots, e_r(0)) \in E$ implies $x(0) \in Z_e$.
 System (\ref{sysappb}) augmented with (\ref{exo1}) can be
 regarded as a particular case of system (\ref{sysappA}) of the Appendix,
 to which Lemma \ref{lemmaAppB}  applies. In particular by property (b) of the
 latter, there is a number $g^\star>0$ such that for all $g \geq g^\star$
 the positive orbit of ${\bf Z} \times Z_e$ under the flow of (\ref{exo1}) --
 (\ref{sysappb}) is bounded. As a
 consequence, the $\omega$-limit set $\tilde {\cal A}_0$ of ${\bf Z} \times
 Z_e$, is a nonempty, compact, invariant set
 which uniformly attracts ${\bf Z} \times Z_e$. We prove now that
 $\tilde {\cal A}_0 = {\cal A}_0 \times
 \{0\}$. To this end note first of all that ${\cal A}_0  \times
 \{0\}$ by construction is contained in ${\tilde A}_0$. Moreover, $x$ is necessarily 0 at any point of ${\cal A}_0$. In fact suppose, by
 contradiction, that there is a point $({\bf z},x)$ of $\tilde {\cal A}_0$
 with $x \neq 0$. As $g>0$ and $A$ is Hurwitz,
 it follows that the trajectory $x(t)$ of (\ref{exo1}) -- (\ref{sysappb})
 originating from $({\bf z},x)$ is unbounded
 in backward time, which contradicts the fact that $\tilde {\cal A}_0$
 is a compact invariant (in particular in backward time) set. Finally, since
 ${\cal A}_0$ is the $\omega$-limit set of ${\bf Z}$ under the
 flow of (\ref{zerodynpla}), we can conclude that necessarily
  $\tilde {\cal A}_0 = {\cal A}_0 \times \{0 \}$. This in particular proves claim
  (ii)'. Claim (iii)', namely exponential attractivity of $\tilde {\cal
A}_0$, is an easy consequence
 of property (a) of Lemma \ref{lemmaAppB} and of the fact that
 the lower subsystem of (\ref{sysappb}) is exponentially
 stable. To prove claim (iv)', note that (\ref{friend1}),
 (\ref{tildeqq}), (\ref{friend2}) imply
 $
 \left . \tilde c(\tilde {\bf z}) \right |_{{\cal Z} \times \{0\}}
 = c({\bf z})\,.
 $
 From this claim (iv)' immediately follows by assumption (iv),
 taking as $\tilde \tau(\tilde {\bf z})$ any differentiable
 function such that
 $
 \left . \tilde \tau(\tilde {\bf z}) \right |_{{\cal A}_0 \times
 \{0\}} = \left . \tau({\bf z}) \right |_{{\cal A}_0}\,.
 $
 This completes the proof. $\triangleleft$\end{proof}

 Motivated by the previous considerations and result, in what follows we focus our
 attention on the case in which $r=1$, i.e. on the special case in
 which system (\ref{pla1}) is a system of the form
 \beeq{\label{plard1}\ba{rcl} \dot z &=& f_0(\varrho,w,z) +
 f_1(\varrho,w, z,e_1)e_1\\
 \dot e_1 &=& q(\varrho,w,z,e_1) + u\\
 e &=& e_1\\y &=& e_1
 \ea
 }
 and we assume that assumptions (i)\,-\,(ii)\,-\,(iii)\,-\,(iv) hold.

\section{ The adaptive internal model}

\subsection{The structure of the regulator}

The proposed regulator is a system of the form
 \beeq{\label{control}\ba{rcl} u &=& \xi_1 + v\\[2mm]
 \dot \xi &=& A\xi + \phi(\xi_1) + \Omega(\xi_1)\hat \theta
 +H(X,\xi_1)v - M(X) \mbox{dzv}_\ell(\hat \theta)\\[2mm]
 \dot {\hat \theta} &=& \beta(X,\xi_1) v -\mbox{dzv}_\ell(\hat \theta)\\[2mm]
 \dot X &=& FX + G\Omega(\xi_1)\ea
 }
 in which $\xi_1$ denotes the first component of $\xi$, the matrix
$X$ is a $(d-1)\times q$ matrix, $M(X)$ is a $d\times q$ matrix
defined as
 \[
 M(X) = \qmx{0 \cr X\cr}\,,
 \]
 while the vectors $H(X,\xi_1),
\beta(X,\xi_1)$ and the matrices $F, G$ have the form described
below.  The function dzv$_\ell(\cdot)$ is defined as
\[
    \mbox{dzv}_\ell(\mbox{col}(s_1,\ldots,s_q)) = \mbox{col}(\mbox{dz}_\ell(s_1), \ldots,
    \mbox{dz}_\ell(s_q))
\]
in which dz$_\ell(\cdot)$ is any continuously differentiable
function satisfying
 \beeq{\label{deftas}
\mbox{dz}_\ell(x) =
 \left \{ \ba{ll} 0 & \quad \mbox{if } |x|
 \leq \ell\\
 x & \quad \mbox{if } |x| \geq \ell +1
  \ea \right .
 }
 and the amplitude $\ell$ of the {\em dead-zone} is chosen so that
\[
\ell > \max_{\varrho \in P} |\theta(\varrho)|\,.
\]
 This controller can be viewed as a ``copy" of
(\ref{adaptim}), corrected by an ``innovation term", augmented
with an ``adaptation law" for ${\hat \theta}$ and with a ``filter"
which generates the ``auxiliary state" $X$. The additional input
$v$, which is a ``stabilizing control", will eventually be taken
as $v = -k y$.

Following the theory of adaptive observers of \cite{BaGe} and
\cite{MaTo92}, the functions $H(X,\xi_1), \beta(X,\xi_1)$ and the
matrices $F, G$ of (\ref{control}) are chosen as follows. Define
new variables \beeq{\label{newvar}\ba{rcl}
\tilde \theta &=& \hat \theta - \theta(\varrho)\\[2mm] \eta &=& \xi -M(X)\tilde
\theta\,. \ea} (note that $\eta_1 = \xi_1$) and observe that, in
the new variables, the second equation of (\ref{control}) reads as
follows (for convenience, we omit the arguments $(X,\xi_1)$ in $H$
and $\beta$ and the argument $X$ in $M$)
 \beeq{\label{secondeq}
 \ba{rcl} \dot \eta &=& A(\eta + M\tilde
 \theta)  + \phi(\xi_1) + \Omega(\xi_1)( \theta(\varrho) + \tilde\theta) +
 Hv - \dot M
 \tilde \theta - M\beta v\\[2mm]
 &=& A\eta + [AM + \Omega(\xi_1) - \dot M]\tilde \theta +
 [H-M\beta]v + \phi(\xi_1) + \Omega(\xi_1) \theta(\varrho)\,.\ea
 }
 The third equation, instead, becomes trivially
 \[
 \dot {\tilde \theta} = \beta v-\mbox{dzv}_\ell\,(\tilde \theta +
 \theta(\varrho))\,.
 \]
 The choices of $H(X,\xi_1), \beta(X,\xi_1)$ and of $F, G$ are
meant to simplify the terms
\[
[AM + \Omega(\xi_1) - \dot M]\tilde \theta + [H-M\beta]v
\] in the expression (\ref{secondeq}).
First of all, note that choosing
\[
H = M\beta + K
\]
with $K$ a constant vector (whose expression will be determined
later), the second term becomes equal to $Kv$. As for the first
term, the idea is to impose that
\[
[AM + \Omega(\xi_1) - \dot M]\tilde \theta= b \beta\tr \tilde
\theta
\]
in which $b$ is a $d\times 1$ fixed vector. The identity in
question holds  if $M$ satisfies
\[
\dot M = (A-bCA)M + (I-bC)\Omega(\xi_1)
\]
and $\beta$ is taken as \[ \beta\tr = CAM + C\Omega(\xi_1)\,.
\]

In this way, the second equation of (\ref{control}) takes the
simplified form \beeq{\label{MarTom} \dot \eta = A\eta + b\beta\tr
\tilde \theta + Kv+\phi(\eta_1) + \Omega(\eta_1)
\theta(\varrho)\,, } on which we will return later. To show that
the required differential equation for $M$ can be enforced, pick a
column vector $b = {\rm col}(1,b_2, \ldots, b_d)$. Then, bearing
in mind the definition of $M$, it is easily realized that the
required differential equation holds if the matrices $F$ and $G$
in the differential equation for $X$ have the form (see
\cite{MaTo92}) \beeq{\label{effegi} F = \qmx{-b_2 & 1 & \cdots & 0
& 0\cr \cdot & \cdot & \cdots & \cdot & \cdot\cr -b_{d-1} & 0 &
\cdots & 0 & 1\cr -b_{d} & 0 & \cdots & 0 & 0\cr}, \qquad G =
\qmx{-b_2 & 1 & \cdots & 0 & 0 & 0\cr \cdot & \cdot & \cdots &
\cdot & \cdot & \cdot\cr -b_{d-1} & 0 & \cdots & 0& 1 & 0\cr
-b_{d} & 0 & \cdots & 0 & 0 & 1\cr}.}

In summary,  the quantities $H(X,\xi_1),\beta(X,\xi_1),F,G$ which
appear in the controller (\ref{control}) are determined as
follows: $F$ and $G$ are the matrices in (\ref{effegi}),
$\beta(X,\xi_1)$ is chosen as \beeq{\label{biggamma}
\beta(X,\xi_1) = [CA\qmx{0 \cr X\cr} + C\Omega(\xi_1)]\tr } and
$H(X,\xi_1)$ is chosen as \beeq{\label{acca} H(X,\xi_1) = \qmx{0
\cr X\cr}[CA\qmx{0 \cr X\cr} + C\Omega(\xi_1)]\tr + K\,.} The
vector $b$, whose entries determine the choice of $F$ and $G$ and
the parameter $K$, which appears in the expression of
$H(X,\xi_1)$, will be chosen later.

The controller thus defined determines a closed loop system which,
in the coordinates indicated above, can be written as (recall that
$e_1=e$)
 \beeq{\label{closedloop}\ba{rcl}
 \dot \varrho &=& 0 \\[2mm]
 \dot w &=& s(\varrho,w) \\[2mm]
 \dot z &=& f_0(\varrho, w, z)+f_1(\varrho,w,z,e)e\\[2mm]
 \dot e &=& q(\varrho,w,z,e)+\eta_1+v\\[2mm]
 \dot \eta &=& A\eta + b\beta\tr \tilde \theta + K v + \phi(\eta_1) + \Omega(\eta_1)\theta(\varrho)
 \\[2mm] \dot {\tilde \theta} &=& \beta v -\mbox{dzv}_\ell\,(\tilde \theta + \theta(\varrho))\\[2mm]
\dot X &=& FX + G\Omega(\eta_1)\,,\ea}
 where $\beta$ is a function of $X$ and $\eta_1$. This system, viewed as a system with input
 $v$ and output $e$, has relative degree 1 and its zero dynamics
 are those of
 \beeq{\label{zerodyn}\ba{rcl}
 \dot \varrho &=& 0 \\[2mm]
 \dot w &=& s(\varrho,w) \\[2mm]
 \dot z &=& f_0(\varrho, w,z) \\[2mm]
 \dot \eta &=& A\eta - K[q(\varrho,w,z,0)+\eta_1]+
 b\beta\tr \tilde \theta + \phi(\eta_1) + \Omega(\eta_1)\theta(\varrho)
 \\[2mm] \dot {\tilde \theta} &=& -\beta [q(\varrho,w,z,0)+\eta_1]
 -\mbox{dzv}_\ell\,(\tilde \theta + \theta(\varrho))
 \\[2mm]
 \dot X &=& FX + G\Omega(\eta_1)\,.
 \ea
 }
 The intuition suggests that  if the latter have convenient asymptotic properties,
 in particular possess
 a locally exponentially stable compact attractor, an additional control of the form $v=-ke$,
 (with large $k>0$) should be able to solve the problem of output regulation.
 Thus, in following subsection, the asymptotic properties of
 (\ref{zerodyn}) will be studied.

\subsection{Trajectories of (\ref{zerodyn}) are bounded}

In studying the asymptotic properties of this system, it is
convenient to take advantage of the ``immersion" assumption (iii)
introduced above. Specifically, suppose that the initial
conditions for $\varrho,w,z$ are taken in the set ${\bf Z} $, a
subset of a set ${\cal Z}$ which by  hypothesis is positively
invariant for the subsystem formed by the top three equations of
(\ref{zerodyn}). Thus, for any  of such initial conditions and for
any $t\ge 0$, the function $\tau(\varrho, w(t),z(t))$ is well
defined and it is legitimate to consider the change of variables
\[
\chi = \eta - \tau(\varrho,w,z)\,.
\]
This transforms system (\ref{zerodyn}) in a system of the form
(use here (\ref{adaptimmers1}) and (\ref{adaptimmers0}) which hold
on ${\cal A}_0 \subset {\cal Z}$ )
 \beeq{\label{zerodyn2}
 \ba{rcl}
 \dot \varrho &=& 0 \\[2mm]
 \dot w &=& s(\varrho,w) \\[2mm]
 \dot z &=&  f_0(\varrho,w,z) \\[2mm]
 \dot \chi &=& (A-KC)\chi_1+ b\beta\tr \tilde \theta +
 \Delta(\chi_1,\tau_1,\theta) + e(\varrho, w,z)\\[2mm]
 \dot {\tilde \theta} &=& -\beta\chi_1 -\mbox{dzv}_\ell\,(\tilde \theta +
 \theta(\varrho))\\[2mm]
 \dot X &=& FX + G\Omega(\chi_1+\tau_1)\,,\ea} in which
 \[
 \Delta(\chi_1,\tau_1,\theta) = \phi(\chi_1+\tau_1) - \phi(\chi_1)
 + [\Omega(\chi_1+\tau_1) - \Omega(\chi_1)]\theta(\varrho)
 \]
 is a term which vanishes at $\chi_1=0$ and
 \beeq{\label{defe}\ba{rcl}
 e(\varrho, w,z) &=& K(c(\varrho, w,z)-\tau_1(\varrho, w,z)) + A \tau(\varrho, w,z)) + \phi(\tau_1(\varrho, w,z))
 + \Omega(\tau_1(\varrho, w,z)) \theta(\varrho) \\[2mm] &&\dst \qquad-
 \,{\partial \tau \over \partial {z}} {
 f}_0(\varrho, w, z)
 - {\partial \tau \over \partial {w}} {
 s}(\varrho, w)\ea
 }
 is a term vanishing on ${\cal A}_0$.
 In particular note that, since $\phi(y)$ and $\Omega(y)$ can be taken to be
globally Lipschitz and $\theta$ ranges over a compact set, there
exists a number $L$ such that
\[
|\Delta(\chi_1,\tau_1,\theta)| \le
L_\phi\,|\chi_1|+L_\Omega|\chi_1|\,|\theta|\le L|\chi_1|\] for all
$\chi_1,\tau_1,\theta$.

The idea is now to choose the $b_i$'s and $K$ so that system
(\ref{zerodyn2}) has certain desirable asymptotic properties. To
this end, let the $b_i$ be such that the polynomial
 \beeq{ \label{polync}
 p(\lambda) = \lambda^{d-1} + b_2\lambda^{d-2} + \cdots +
 b_{d-1}\lambda + b_d
 }
 has $d-1$ distinct roots with negative real part. As a consequence
 the matrix $F$ in the bottom equation of
(\ref{zerodyn2}) is Hurwitz (and has distinct eigenvalues). This,
in view of the assumptions on the top three equations, suggests
that the asymptotic properties of (\ref{zerodyn2}) are entirely
determined by those of the fourth and fifth equation.

As indicated in \cite[Theorem 2.1]{MaTo92}, the appropriate choice
for $K$ in (\ref{MarTom}) is \beeq{\label{kappa} K = Ab + \lambda
b } in which $\lambda >0$. To see why this is the case note first
of all that, using a little algebra, it is not difficult to prove
the following.

\begin{lemma}\label{MaTo} Choose $K$ as in (\ref{kappa}) and set
\[
T = \qmx{1 & 0 \cr \hat b & I\cr}, \qquad \hat b = - \qmx{b_2 \cr
\cdot \cr b_{d-1} \cr b_{d} \cr}.
\]
Then
\[
T(A-KC)T^{-1} = \qmx{-\lambda & \hat c \cr 0 & F\cr}, \quad Tb
=\qmx{1 \cr 0\cr}, \quad CT^{-1} = C\,,
\]
in which $\hat c = \qmx{1 & 0 & \cdots & 0\cr}$ and $F$ is the
matrix defined in (\ref{effegi}).
\end{lemma}

>From this fact, standard arguments can be invoked to claim
boundedness of the trajectories of (\ref{zerodyn2}). In fact, the
following result holds.

\begin{lemma}\label{LM1} Suppose assumptions (i), (ii), (iv) hold. There is a number $\lambda^\ast$ such
that, if $\lambda \ge \lambda^\ast$, all trajectories of
(\ref{zerodyn2}) are bounded.
\end{lemma}

\begin{proof}
First of all, recall that, by assumption (ii), $(\varrho, w(t),
z(t)) \in {\cal Z}$ for all $t \geq 0$, where ${\cal Z}$ is a
compact set. Thus, looking at the expression (\ref{defe}) of
$e(\varrho,w,z)$, it is seen that there exist a number $\bar e$
(depending on the design parameter $\lambda$ because the latter
appears in $K$) such that \beeq{ \label{ebar}
 |e(\varrho, w(t), z(t))| \leq \bar e \qquad \forall \; t \geq 0\,.
 }
Observe also that, so long as trajectories of (\ref{zerodyn2})
exist on some interval $[0,T]$, $|X(t)|$ is bounded, by a number
which only depends on $|X(0)|$ (because $|\Omega(\cdot)|$, having
compact support, is bounded by some fixed number). As a
consequence, also $|\beta(t)|$ is bounded, again by a number which
only depends on $|X(0)|$.  Thus, system (\ref{zerodyn2}) cannot
have finite escape times.

This being the case, to prove the Lemma it remains to show that
also $\chi$ and $\tilde \theta$ are bounded. To this end, let
$\chi$ be partitioned as $\chi = {\rm col}(\chi_1,\chi_2)$, in
which $\chi_2$ is a $(d-1)\times 1$ vector and change $\chi_2$
into
 \[
 \zeta = \hat b \chi_1 + \chi_2\,.
 \]
 In this way, the fourth and fifth equations of system
 (\ref{zerodyn2}) are changed into
 \beeq{\label{MarTom2}\ba{rcl}
 \dot \chi_1 &=& -\lambda \chi_1 + \hat c\zeta +
 \beta\tr \tilde \theta + C\Delta(\chi_1,\tau_1,\theta) + C e(\varrho, w, z)\\[2mm]
 \dot \zeta &=& F\zeta + \qmx{\hat b & I \cr}\Delta(\chi_1,\tau_1,\theta)
 +\qmx{\hat b & I \cr} C e(\varrho, w, z)\\[2mm]
 \dot {\tilde \theta} &=& -\beta \chi_1 - \mbox{dzv}_\ell(\tilde \theta + \theta(\varrho))\,.
 \ea}

 With this in  mind, choose for (\ref{MarTom2}) the Lyapunov function
 \beeq{\label{Vdef}
 V(\chi_1,\zeta,\tilde \theta) = \chi_1^2 + \zeta\tr P \zeta +
 \tilde \theta\tr \tilde \theta\,,
 }
 in which $P$ is the positive definite solution of $PF+F\tr P = -
 I$. This yields
 \beeq{\label{dotv.1}\ba{rcl} \dot V &=& -2\lambda
 \chi_1^2 + 2 \chi_1\hat c\zeta + 2 \chi_1\beta\tr \tilde \theta +
 2\chi_1 C\Delta(\chi_1,\tau_1,\theta) + 2\chi_1 C e(\varrho, w, z)\\[2mm]
 && - |\zeta|^2 + 2 \zeta\tr P \qmx{\hat b &
 I\cr}\Delta(\chi_1,\tau_1,\theta) + 2 \zeta\tr P \qmx{\hat b &
 I\cr} C e(\varrho, w, z)\\[2mm]
 && - 2\tilde \theta\tr \beta\chi_1 - 2 \tilde \theta\tr
 \mbox{dzv}_\ell(\tilde \theta + \theta(\varrho)) \\[2mm]
 &\le& -2\lambda \chi_1^2 - |\zeta|^2 - 2 \tilde \theta\tr
 \mbox{dzv}_\ell(\tilde \theta + \theta(\varrho)) + L_1|\chi_1|^2 +
 L_2|\chi_1|\,|\zeta|\\[2mm]
 && +L_3|\chi_1||e(\varrho, w, z)| + L_4 |\zeta||e(\varrho, w,
 z)|
 \ea}
 in which $L_i$, $i=1\ldots,4$ are suitable positive constants.
 By completing the squares and using (\ref{ebar}), we obtain
 \beeq{\label{Vdot1}
 \dot V \le -(2\lambda -L_1 + {\dst 1 \over \dst 2}L_2^2) \chi_1^2  - {\dst 1 \over 2}|\zeta|^2
 - 2 \tilde \theta\tr \mbox{dzv}_\ell(\tilde \theta + \theta(\varrho))
 +L_3|\chi_1|\bar e + L_4 |\zeta|\bar e\,.
 }
Bearing in mind the definition (\ref{deftas}) and the choice of
$\ell$, observe that \beeq{\label{tas2}
  \tilde \theta\tr \mbox{dzv}_\ell(\tilde \theta +
 \theta(\varrho)) \geq 0 \qquad \mbox{for all } \tilde \theta \in
 \Real^q\quad \mbox{and } \quad \varrho \in P\,.
 }
 It is also easy to check that for any $\delta >  \sqrt{q}(2 \ell+1)$ there is a positive
 number $c_1$ such that
 \beeq{\label{tas1}
 |\tilde \theta| \geq \delta \qquad \Rightarrow \qquad
 2 \tilde \theta\tr \mbox{dzv}_\ell(\tilde \theta + \theta(\varrho)) \geq
 c_1 |\tilde \theta|^2 \qquad \mbox{for all } \tilde \theta \in
 \Real^q\quad \mbox{and } \quad \varrho \in P\,.
 }
Pick a value of $\lambda$ large enough so that $\bar \lambda:=2
\lambda - L_1 - L_2^2/2>0$. Inequality (\ref{Vdot1}), in view of
property (\ref{tas1}), yields
 \[
 |\tilde \theta| \geq \delta \qquad \Rightarrow \qquad \dot V \leq
 -c_2  |(\chi_1, \zeta, \tilde \theta)|^2  + c_3 |(\chi_1, \zeta, \tilde
 \theta)|\bar e
 \]
 in which $c_2 = \min \{\bar \lambda, \dst {1\over 2}, c_1\}$ and $c_3 = 3 \max \{L_3,
 L_4\}$. From this, it is seen that
 \beeq{\label{Re1}
  |\tilde \theta| \geq \delta \quad \mbox{and} \quad
  |(\chi_1, \zeta, \tilde \theta)| > {c_3 \over c_2} \, \bar e
 \qquad \Rightarrow \qquad \dot V <0\,.
 }
Property (\ref{tas2}), on the other hand, yields
 \[
 \dot V \leq - c_2 |(\chi_1, \zeta)|^2 + c_3 |(\chi_1, \zeta)| \,
 \bar e
 \]
 from which it is seen that
 \beeq{\label{Re2}
 |(\chi_1, \zeta)| > {c_3 \over c_2} \, \bar e
 \qquad \Rightarrow \qquad \dot V <0\,.
 }

 We show now that a combination of (\ref{Re1}) and (\ref{Re2})
 yields the desired result, namely the boundedness of
 $(\chi_1(t), \zeta(t), \tilde \theta (t))$. As a matter of fact set
 \[
 r := \sqrt{\delta^2 + \Bigl({c_3 \over c_2} \,\bar e\Bigr)^2}
 \]
and note that, since
 \[
 |(\chi_1, \zeta, \tilde \theta)| > r \qquad \Rightarrow \qquad |(\chi_1, \zeta, \tilde \theta)|
 > {c_3 \over c_2} \bar e
 \]
 and
 \[
 |(\chi_1, \zeta, \tilde \theta)| > r \qquad \Rightarrow \qquad
 |(\chi_1, \zeta)| > {c_3 \over c_2} \bar e \quad \mbox{or} \quad
 |\tilde \theta| >\delta\,,
 \]
 relations (\ref{Re2}) and (\ref{Re1}) imply
 \[
 |(\chi_1, \zeta, \tilde \theta)| > r \qquad \Rightarrow \qquad
 \dot V <0\,.
 \]
 From this, bearing in mind the fact that $V(\chi_1,\zeta,\tilde \theta)$ is a quadratic form,
  the result follows by standard arguments.
 $\triangleleft$
\end{proof}

 We can therefore draw the following conclusion about system (\ref{zerodyn}). Let the initial
conditions $\eta(0),\tilde \theta (0), X(0)$  be taken in fixed
compact sets ${\bf H},\, \Theta,\, {\bf X}$.\,\footnote{Recall
that $\tilde \theta(t) = \hat \theta(t) - \theta(\varrho)$ and $
\eta(t) = \xi(t) - \beta(X(t),\xi_1(t))\tilde \theta(t)$. Thus, to
establish boundedness of trajectories when $\xi (0)$, $\hat
\theta(0)$ and $X(0)$ are taken in fixed compact sets it suffices
to consider the case in which $\tilde \theta(0)$, $\eta(0)$ and
$X(0)$ are taken in fixed compact sets.} Then, the positive orbit
of the set
\[
{\bf B} = {\bf Z} \times {\bf H} \times \Theta \times {\bf X}
\]
under the flow of (\ref{zerodyn}) is bounded. As a consequence
$\omega({\bf B})$, the $\omega$-limit set of ${\bf B}$ under the
flow of (\ref{zerodyn}), is a non-empty, compact and invariant
set, which uniformly attracts  all trajectories of (\ref{zerodyn})
with initial conditions in ${\bf B}$.

\subsection{The limit set of (\ref{zerodyn})}\label{sec3}

We proceed now to investigate the structure of the set
$\omega({\bf B})$. To this end, we look at the equivalent system
(\ref{zerodyn2}), we note that the three top equations are
independent of the bottom ones and we rewrite them in compact form
as in (\ref{zerobold}) (and consistently we rewrite the term
$e(\varrho, w,z)$ as $e({\bf z})$ and $\tau(z,w,\varrho)$ as
$\tau({\bf z})$). In particular, because of the special triangular
structure of (\ref{zerodyn2}), we note that if $({\bf
z},\chi,\tilde \theta, X)$ is a point of $\omega({\bf B})$,
necessarily ${\bf z}$ is a point in the $\omega$-limit set of
${\bf Z}$ under the flow of (\ref{zerobold}), that is, ${\bf z}$
is a point of ${\cal A}_0$. This implies that on $\omega({\bf B})$
we have $e({\bf z})=0$ and thus system (\ref{zerodyn2}) simplifies
as
 \beeq{\label{zerodyn.2}
 \ba{rcl}\dot
 {\bf z} &=& {\bf f}_0({\bf z})\\[2mm]\dot \chi &=& (A-KC)\chi + b\beta\tr \tilde
 \theta +
 \Delta(\chi_1,\tau_1,\theta) \\[2mm] \dot {\tilde \theta} &=& -\beta
 \chi_1 - \mbox{dzv}_\ell(\tilde \theta + \theta(\varrho))\\[2mm]
 \dot X &=& FX + G\Omega(\chi_1+\tau_1({\bf z}))\,.\ea
 }

 What we will be able to prove in the following is that
 on points of $\omega({\bf B})$ necessarily $\chi=0$, $\tilde \theta=0$ and
 the value of $X$ is entirely determined by the properties of the
 system
 \beeq{\label{zerodyn3} \ba{rcl} \dot
 {\bf z} &=& {\bf f}_0({\bf z})\\[2mm]
 \dot X &=& FX + G\Omega(\tau_1({\bf z})),\ea
 }
 in which $\tau_1({\bf z})$ is the obvious abbreviated notation for
$\tau_1(z,w,\varrho)$.
 To this end, though, an extra hypothesis is needed, which will be explained
 after having shown an interesting feature of the system in
 question.

\begin{lemma}\label{ssX}
The graph of the map
\[
\ba{rcccl}
 \sigma &:& {\cal  A}_0 &\to& \Real^{(d-1)\times q}\\
 && {\bf z} & \mapsto & \dst \int_{-\infty}^0
 e^{-Fs}G\Omega(\tau_1({\bf z}(s,{\bf z}))) ds
 \ea
\]
 is invariant for (\ref{zerodyn3}).
\end{lemma}

\begin{proof}
 Let ${\bf z}(t,{\bf z}_0)$ denote the solution of (\ref{zerobold}) passing
through ${\bf z}_0$ at time $t=0$ and note that, if ${\bf z}_0 \in
{\cal A}_0$, then ${\bf z}(t,{\bf z}_0) \in {\cal A}_0$ for all
$t$ (thus, in particular, since ${\cal A}_0$ is compact, $|{\bf
z}(t,{\bf z}_0)|$ is bounded by a number which depends only on
${\cal A}_0$). Then, since $F$ is a Hurwitz matrix, the map
$\sigma(\cdot)$ is well defined. As simple calculation shows that
\[
 \sigma({\bf z}(t,{\bf z}_0)) =e^{F t} \sigma({\bf z}_0) +
 \int_0^t e^{F(t- s)} G \Omega(\tau_1({\bf z}(s,{\bf z}_0)))
 ds\,.
\]
This shows that
 \[
 {\rm graph}(\sigma ) =\{({\bf z},X): {\bf z}\in {\cal A}_0,
 X=\sigma({\bf z})\}
 \]
is invariant for (\ref{zerodyn3}). $\triangleleft$ \end{proof}

\begin{remark}
Consider the restriction of (\ref{zerodyn3}) to ${\cal A}_0 \times
\Real^{(d-1) \times q}$. Since  the graph of $\sigma(\cdot)$ is
invariant for (\ref{zerodyn3}), changing $X$ into $\tilde X= X -
\sigma({\bf z})$, yields \[ \ba{rcl} \dot
{\bf z} &=& {\bf f}_0({\bf z})\\[2mm]
\dot {\tilde X} &=& F{\tilde X}\,.\ea\] We see from this that the
solution $X(t)$ of (\ref{zerodyn3}) passing through $X_0$ at time
$t=0$ can be expressed as
 \beeq{ \label{ssnl} X(t) = e^{Ft}[X_0- \sigma({\bf
z}_0)] + \sigma({\bf z}(t,{\bf z}_0))\,. \qquad \triangleleft }
\end{remark}

We introduce now an additional hypothesis, reminiscent of the
classical hypothesis of {\em persistence of excitation}.

\bigskip\noindent
{\em Assumption (v)}\,: Consider the map $\varphi: {\cal A}_0 \to
\Real^{q\times 1}$ defined as \[ \varphi:  {\bf z}
\;\;\;\mapsto\;\;\; \beta(\sigma({\bf z}),\tau_1({\bf z}))\] It is
assumed that for any initial condition ${\bf z}_0 \in {\cal A}_0$
 the identity
\[
\gamma\tr \varphi({\bf z}(t,{\bf z}_0)) = 0, \qquad \mbox{for all
$t\in \Real$}\] implies $\gamma=0$. $\triangleleft$

\begin{remark} In other words, the assumption of ``persistency of
excitation", in the present context, is spelled as follows: for
any initial condition ${\bf z}_0 \in {\cal A}_0$, the $q$ outputs
of the autonomous system
\[\ba{rcl}
\dot {\bf z} &=& {\bf f}_0({\bf z})\\[1mm]
\varphi &=& \beta(\sigma({\bf z}),\tau_1({\bf z}))\ea\] are
linearly independent functions, on the entire time axis.
$\triangleleft$ \end{remark}

\bigskip
Under this hypothesis, the set $\omega({\bf B})$ assumes a very
simple structure. As a matter of fact, the following result holds.

\begin{lemma} \label{Lemmaomegaset}Suppose that, in addition to
assumptions (i), (ii), (iv), also assumption (v) holds. Then the
values of $\chi$ and $\tilde \theta$ on any point of $\omega({\bf
B})$ are necessarily zero.\end{lemma}

\begin{proof}  By contradiction, suppose
a point ${\bf p}=({\bf z},\chi_0,\tilde \theta_0, X)$ with either
$\chi_0\ne 0$ or $\tilde \theta_0\ne 0$ is in $\omega({\bf B})$.
Since $\omega({\bf B})$ is compact and invariant, in particular in
backward time, the backward trajectory of (\ref{zerodyn.2})
starting at this point is bounded. Along this trajectory, the
function
 \[
 V(t) := V(\chi_1(t),\zeta(t),\tilde \theta(t))
 \]
 in (\ref{Vdef}) satisfies $V(t) \le C$ for all $t\le 0$, for some $C>0$.
 Moreover, since $e({\bf z})=0$ on $\omega({\bf B})$, the same
 computations indicated in the proof of Lemma \ref{LM1} show that

 \[
 \dot V(t) \leq -(2 \lambda - L_1 - {1 \over 2} L_2^2)|\chi_1(t)|^2
 -{1 \over 2} |\zeta(t)|^2 - 2\tilde \theta(t) \mbox{dzv}_{\ell}(\tilde \theta(t) + \theta(\varrho))
 \]
 in which $L_1$ and $L_2$ are the same constants introduced in the proof of
 Lemma \ref{LM1}.  From this, using property (\ref{tas2}),
  it turns out that if $\lambda \geq \lambda^\star$ (where
  $\lambda^\star$  is the same as in Lemma \ref{LM1}) then $V(t)$ is
 non-increasing along trajectories.  As consequence, since $V(t)$ is bounded, that there
 must exist a finite number $V_\alpha$ such that
 \[\lim_{t\to -\infty}V(t)=V_\alpha\,.\]
The trajectory in question is attracted, in backward time, by its
own $\alpha$-limit set $\alpha({\bf p})$, which, as it is well
known, is nonempty, compact and invariant. Moreover, by
definition, the function $V(\chi_1,\zeta,\tilde \theta)$ has the
same value $V_\alpha$ at any point of $\alpha({\bf p})$.

Now, as in the classical proof of LaSalle's invariance principle,
pick an initial condition $\hat {\bf p}$ in the set $\alpha( {\bf
p})$ and consider the corresponding trajectory of
(\ref{zerodyn.2}), which remains in $\alpha( {\bf p})$  for all
times. Along such trajectory, $V(t)$ is constantly equal to
$V_\alpha$ and hence
\[
\chi_1(t) = 0, \qquad \zeta(t)=0, \qquad
 \mbox{dzv}_\ell(\tilde \theta(t) + \theta(\varrho))=0\qquad
\mbox{for all $t\in \Real$.}\] Entering these constraints in
(\ref{zerodyn.2}), and observing that the vector $b$ is nonzero,
it is seen that necessarily
\[\ba{l}
\tilde \theta \tr \beta  =0\\[2mm]
\dot{\tilde \theta} = 0\\[2mm]
\dot X = FX + G\Omega(\tau_1({\bf z}))\,.\ea\] The second
condition shows that $\tilde \theta(t)$ is a constant, say $\tilde
\theta^\ast$, along such trajectory. The third condition, says
that $X(t)$ is a solution of
\[
\dot X = FX + G\Omega(\tau_1({\bf z}))
\]
Now, since $F$ is Hurwitz and has distinct eigenvalues (because so
are the roots of the polynomial (\ref{polync})), it is seen from
(\ref{ssnl}) that $X(t)$ is bounded for $t\le 0$ only if $X(0) =
\sigma({\bf z}(0))$, where $\sigma(\cdot)$ is the map introduced
in Lemma \ref{ssX}, in which case $X(t) = \sigma({\bf z}(t))$.
Since $X(t)$ has to be bounded because $\alpha({\bf p})$ is
compact, it follows that $X(t)$ is necessarily equal to
$\sigma({\bf z}(t))$. This being the case,
 bearing in mind  the expression of $\beta$ and the definition of the
map $\varphi(\cdot)$, the first condition shows that necessarily
\[
 ( \tilde\theta^\ast)\tr \varphi({\bf z}(t))= 0, \qquad \mbox{for all $t\in
\Real$}.\] Thus, in view of the assumption of persistency of
excitation, it follows that $\tilde \theta^\ast =0$. It is seen in
this way that $(\chi_1,\zeta,\tilde \theta)=(0,0,0)$ at any point
of $\alpha( {\bf p})$, and this  proves that $V_\alpha =0$. But
this is a contradiction, because $V(t)$ is non-increasing along
trajectories and $V(0)$ is strictly positive, if either $\chi_0
\ne 0$ or $\tilde \theta_0\ne 0$. $\triangleleft$
\end{proof}

To complete the analysis, it remains to determine the values of
$X$ on points of $\omega({\bf B})$. Knowing that $\chi_1=0$ on any
of such points, it follows from the previous analysis and in
particular from Lemma \ref{ssX} that $X=\sigma({\bf z})$.
Altogether, bearing in mind how system (\ref{zerodyn}) and system
(\ref{zerodyn.2}) are related, the following conclusion holds.

\begin{proposition}\label{PR1} Under the assumptions (i),(ii),(iv) and (v)
the set $\omega({\bf B})$ is the graph of a
continuous map defined on ${\cal A}_0$. Any point of $\omega({\bf
B})$ is a point $({\bf z},\eta,\tilde \theta,X)$ in which ${\bf
z}\in {\cal A}_0$ and
\[\eta=\tau({\bf z}), \qquad
\tilde \theta = 0, \qquad X = \sigma({\bf z})\,.\]
\end{proposition}

\subsection{Exponential attractivity of the limit set of (\ref{zerodyn})}\label{Sec44}

Finally, we prove that the set $\omega({\bf B})$ is also locally
exponentially attractive for the trajectories of the zero dynamics
(\ref{zerodyn}) of system (\ref{closedloop}), if so is the set
${\cal A}_0$ for the trajectories of (\ref{zerobold}). This fact
is formalized in the next proposition.
%

\begin{proposition}\label{PR2} Suppose that, in addition to assumptions (i)-(ii)-(iv)
and (v), also assumption (iii) holds.  Then $\omega({\bf B})$ is
locally exponentially attractive for
(\ref{zerodyn}).\end{proposition}

\begin{proof}
  Consider again the equivalent system (\ref{zerodyn2}), let the compact notation
  $\dot{\bf z} = {\bf f}_0({\bf z})$
  be used for the
  first three equations and let the variables $(\chi_1, \zeta)$,
  introduced in the proof of Lemma \ref{LM1}, replace $\chi$.
  Let $\bar \sigma : {\cal Z} \;\; \to
  \;\; \Real^{(d-1)\times q}$ be any continuously differentiable map
  which agrees on ${\cal A}_0$ with the map
  $\sigma$ introduced in Lemma \ref{ssX}, and change $X$ into
  $\tilde X = X - \sigma({\bf z})$. In this way,   the last
  equation of (\ref{zerodyn2}) is transformed into an equation of the form
  \[
  \dot {\tilde X} = F \tilde X + Q({\bf z}) + R(\chi_1, {\bf z})
  \]
 in which
 \[
 Q({\bf z})= F \bar \sigma({\bf z}) + G \Omega(\tau_1({\bf z})) -
 {\partial \bar \sigma \over \partial {\bf z}} {\bf f}_0({\bf z})
 \]
 is a (matrix-valued) function vanishing on ${\cal A}_0$ while
 \[
 R(\chi_1, {\bf z}) = G[\Omega(\chi_1 + \tau_1({\bf z})) - \Omega(\chi_1)]
 \]
 is vanishing for $\chi_1=0$ for all ${\bf z} \in {\cal Z}$. Let
 $\tilde X_i, Q_i({\bf z}), R_i(\chi_1, {\bf z})$ denote the
 $i$-th columns of $\tilde X, Q({\bf z}), R(\chi_1, {\bf z})$.
 Setting ${\bf x} = \mbox{col}(\chi_1, \zeta, \tilde \theta, \tilde
 X_1, \ldots,\tilde X_{q})$, system (\ref{zerodyn2}) can be conveniently rewritten as
\[\ba{rcl}
\dot {\bf z} &=& {\bf f}_0 ({\bf z})\\
  \dot {\bf x} &=& {\bf g}({\bf z},{\bf x}) + \nu({\bf z})
 \ea\]
in which
\[
 {\nu}({\bf z}) = \left (
 \ba{c}
  C e({\bf z})\\
  \left (\ba{cc} \hat b & I  \ea \right ) e({\bf z})\\
  0\\
  Q_1({\bf z})\\ \vdots \\ Q_q({\bf z})
  \ea \right )
\]
is a vector of functions vanishing on ${\cal
 A}_0$.
Observing that ${\bf g}({\bf z},0) = 0$, set
 \[ {\bf A}({\bf z})=
 {\partial {\bf g}\over \partial {\bf x}}({\bf z},0)
 \]
 and consider the expansion
\[
{\bf g}({\bf z},{\bf x}) = {\bf A}({\bf z}){\bf x} + {\bf h}({\bf
z},{\bf x})\,.\] The matrix ${\bf A}({\bf z})$ is the matrix
\[
 {\bf A}({\bf z}) =
 \left ( \ba{cccccc}
  - \lambda + r_1({\bf z}) & \hat c & \beta\tr(\sigma({\bf z}), \tau_1({\bf z})) & 0& \cdots &0\\
  r_2({\bf z}) & F & 0 & 0& \cdots &0\\
  -\beta(\sigma({\bf z}), \tau_1({\bf z})) & 0 & 0 & 0& \cdots &0\\
  r_{31}({\bf z}) & 0 & 0 & F& \cdots &0\\\cdot & \cdot & \cdot & \cdot& \cdots &\cdot\\r_{3q}({\bf z}) & 0 & 0 & 0& \cdots
  &F
  \ea \right )
 \]
 in which
 \[
 r_1({\bf z}) =C \left[{\partial \Delta
 \over \partial \chi_1 } \right ]_{\chi_1=0}
 \qquad
 r_2({\bf z}) = \left (\ba{cc} \hat b & I  \ea \right )\left[{\partial \Delta
 \over \partial \chi_1 } \right]_{\chi_1=0}
 \qquad
 r_{3i}({\bf z}) = \left[{\partial R_i
 \over \partial \chi_1 } \right]_{\chi_1=0}.
 \]
 Moreover, by construction, the vector ${\bf h}({\bf
z},{\bf x})$ is such that \[ \lim_{|{\bf x}|\to 0}{|{\bf h}({\bf
z},{\bf x})|\over |{\bf x}|}=0,\] uniformly in ${\bf z}$ (as the
latter ranges over a compact set). In this way, system
(\ref{zerodyn2}) is rewritten as \beeq{\ba{rcl} \label{zxsubs}
\dot {\bf z} &=& {\bf f}_0({\bf z})\\ \dot {\bf x} &=& {\bf
A}({\bf z}){\bf x} + {\bf h}({\bf z},{\bf x})+\nu({\bf z})\,.\ea}

 With this in mind, consider now the auxiliary system
 \beeq{\label{tis}
 \ba{rcl}
  \dot {\bf z} &=& {\bf f}_0 ({\bf z})\\
  \dot {\bf y} &=& {\bf A}({\bf z}) {\bf y}
 \ea
  }
 with initial conditions $({\bf z}(0),{\bf y}(0))$ in the compact set ${\bf
 Z}\times {\bf Y}$ where ${\bf Y}=\{{\bf y}:|{\bf y}|\le c\}$, with $c>1$.
Arguments identical to those used in the proof of Lemma
 \ref{LM1} and Lemma \ref{Lemmaomegaset} make it possible to claim
 the existence of a $\lambda^\star>0$ such that for all $\lambda >
 \lambda^\star$ the trajectories $({\bf z}(t), {\bf y}(t))$ are bounded in positive time
 and that \[
 \omega({\bf
 Z}\times {\bf Y}) ={\cal A}_0 \times
 \{0\}\,.\]
  As a matter of fact note that, by assumption, the trajectories
 ${\bf z}(t)$ are such that ${\bf z}(t) \in {\cal Z}$ for all $t \geq 0$.
 As far as the trajectories ${\bf y}(t)$ are concerned,
 consider the candidate Lyapunov function
 \[
 U({\bf y}) = V(\chi_1, \zeta, \tilde \theta) + \sum_{i=1}^q\tilde
 X_i^T P \tilde X_i
 \]
 where $V(\chi_1, \zeta, \tilde \theta)$ is the function defined in (\ref{Vdef}) and
 $P$ is the positive definite solution of $PF + F\tr P = -I$. The
  time derivative of $U({\bf y}(t))$ along the solutions of (\ref{tis})
 can be estimated as
 \[
 \ba{rcl}
 \dot U &=& -\, 2(\lambda - r_1({\bf z})) \chi_1^2 + 2 \hat c \chi_1
 \zeta + 2 \chi_1 \beta \tr(\sigma({\bf z}), \tau_1({\bf z})) \tilde \theta
 - |\zeta|^2 + 2 \zeta\tr P r_2({\bf z}) \chi_1\\[2mm]
 && -\,2 \tilde \theta\tr \beta (\sigma({\bf z}), \tau_1({\bf
 z}))\chi_1 - \sum_{i=1}^q|\tilde X_i|^2 + \sum_{i=1}^q2 \tilde X_i\tr P r_{3i}({\bf z}) \chi_1\\[2mm]
 &\leq&
  -\, 2 (\lambda - \bar r_1) \chi_1^2 - |\zeta|^2 -\sum_{i=1}^q|\tilde X_i|^2 \\[2mm]&&+\, 2 \,\hat
 c\, |\chi_1|\,|\zeta| + 2 \,\bar r_2 \, |P| \,|\zeta| \,|\chi_1| + \sum_{i=1}^q2 \,\bar
 r_3 \,|P| \,|\tilde X_i|\, |\chi_1|
  \ea
 \]
in which $\bar r_i=\max_{{\bf z} \in {\cal Z}} |r_i({\bf z})|$,
 $i=1,2$ and $\bar r_3=\max_{{\bf z} \in {\cal Z}, i=1,\ldots,q} |r_{3i}({\bf z})|$.  Standard arguments can be used to show that
 a large value of $\lambda$ renders $\dot U$
 non positive, from which boundedness of
 ${\bf y}(t)$ follows. Moreover the same arguments of the proof of Lemma
 \ref{Lemmaomegaset} can be repeated to show that, under the condition of
 persistence of excitation expressed by Assumption (v),  points on
 $\omega({\bf Z}\times {\bf Y})$ of (\ref{tis})
 are necessarily characterized by ${\bf y}=0$, from which it follows that
 $\omega({\bf Z} \times {\bf Y})= {\cal A}_0 \times \{ 0 \}$\,.

 We show now that ${\cal A}_0 \times \{ 0 \}$
 is locally exponentially attractive for (\ref{tis}). To this end,
let ${\bf z}(t,{\bf z}_0)$ and ${\bf y}(t,{\bf y}_0)$ denote the
solution pair of (\ref{tis}) passing through ${\bf z}_0$ and,
respectively, ${\bf y}_0$ at time $t=0$. Recall (see section
\ref{basicb}) that ${\cal A}_0 \times \{ 0 \}$ attracts
 the set ${\bf Z} \times
 {\bf Y}$ {\em uniformly}. Therefore, since
 $|{\bf y}| \le {\rm dist}(({\bf z}, {\bf y}),{\cal A}_0 \times \{0\})$,
  for any $\varepsilon>0$ there exists $T_\varepsilon>0$ such
 that
 \beeq{\label{uniformconv}
 |{\bf y}(t,{\bf y}_0)| \leq
 \varepsilon
 \qquad \mbox{for all } \; t \geq T_\varepsilon \quad \mbox{and all}
 \quad
 ({\bf z}_0, {\bf y}_0) \in {\bf Z} \times {\bf Y}.
 }

 With this in mind, let $\delta$ be such that ${\rm dist}({\bf z}_0,{\cal
 A}_0) \le \delta$ implies ${\bf z}(t,{\bf z}_0) \in {\bf Z}$ for all $t\ge
 0$, which is always possible, since ${\cal A}_0$ is stable in the
 sense of Lyapunov for the upper equation of (\ref{tis}). Pick any ${\bf
 z}_0$ within a $\delta$-distance from ${\cal A}_0$ and
regard the bottom equation of (\ref{tis}) as a {\em time-varying
linear} system
 \beeq{\label{ltvs}
 \dot {\bf y} = {\bf A}({\bf z}(t,{\bf z}_0)){\bf y}\,.
 }
 Pick a pair $t\ge t_0 \ge 0$ and let ${\Phi}(t,t_0,{\bf z}_0)$ denote the associated state
 transition matrix (which, of course, depends on the pick of ${\bf z}_0$). By construction, the $i$-th column
 $\phi_i(t,t_0,{\bf z}_0)$ of ${\Phi}(t,t_0,{\bf z}_0)$ is the
 solution of (\ref{ltvs}) which satisfies $\phi_i(t_0,t_0,{\bf z}_0)={\bf
 v}_i$, where ${\bf
 v}_i$ is a vector in which all entries are zero but the $i$-th
 one, which is equal to 1. Consider now again (\ref{tis}) with
 initial conditions ${\bf z}(0) ={\bf z}(t_0,{\bf z}_0)$ and ${\bf
 y}(0)={\bf v}_i$ (note that $({\bf z}(t_0,{\bf z}_0),{\bf v}_i) \in {\bf Z}\times {\bf Y}$). Since (\ref{tis}) is time invariant, we
 observe that ${\bf y}(t,{\bf v}_i)=\phi_i(t+t_0,t_0,{\bf z}_0)$ for all $t\ge 0$.
 Thus, by appealing to (\ref{uniformconv}), it is deduced that, for
 any $\varepsilon$, there exists $T_\varepsilon$ such that
 \[
 |\phi_i(t+t_0,t_0,{\bf z}_0)|\le \varepsilon
\quad \mbox{for all } \; t \geq T_\varepsilon \] and all ${\bf
z}_0$, so long as that ${\rm dist}({\bf z}_0,{\cal
 A}_0) \le \delta$.
This, in turn, by standard results (see e.g. \cite[page 92]{Rugh},
implies the existence of positive numbers $M$ and $a$ (independent
of ${\bf z}_0$) such that
 \beeq{\label{normphi}
 |{\Phi}(t,t_0,{\bf z}_0)| \leq M e^{-a (t - t_0)} \qquad \mbox{for all } \; t \geq t_0 \ge
 0\,,
 }
and all ${\bf z}_0$, so long as that ${\rm dist}({\bf z}_0,{\cal
 A}_0) \le
\delta$.

 By a classical converse Lyapunov theorem (see Theorem 3.12 in
 \cite{Khalil}), we deduce from (\ref{normphi}) the existence of a continuously
 differentiable and symmetric function $\bar P(t)$, of a continuous
 and symmetric function $Q(t)$ and of constants $c_1$, $c_2$ and $c_3$ such that
 \[
 {d \bar P(t) \over dt} + \bar P(t) {\bf A}({\bf z}(t,{\bf z}_0)) + {\bf A}\tr({\bf z}(t,{\bf z}_0)) \bar P (t) = - Q(t)
 \]
 with
 \[
  0 < c_1 I \leq \bar P(t) \leq c_2 I \qquad \mbox{and}
 \qquad Q(t) \geq c_3 I >0
 \]
 for all $t \geq 0$.

 Bearing in mind this result, we return now to the lower
 subsystem of
 (\ref{zxsubs}) which can be more conveniently seen as a {\em time-varying nonlinear} system
 \beeq{
 \dot {\bf x} = {\bf A}({\bf z}(t,{\bf z}_0)) {\bf x} + {\bf h}({\bf z}(t,{\bf z}_0),{\bf x}) + \nu({\bf z}(t,{\bf
 z}_0))\,.
 \label{xssubs}
 }
 In particular note that, as far as the  term ${\bf h}({\bf z}(t,{\bf z}_0),{\bf x})$
 is concerned, for any $\epsilon>0$ there is
 $\delta_\epsilon>0$ such that
 \[
 |{\bf x}| \leq \delta_\epsilon \qquad \Rightarrow \qquad |{\bf h}({\bf z}(t,{\bf z}_0),{\bf x})| \leq  \epsilon |{\bf x}|
 \]
 for all $t \geq 0$ and all ${\bf z}_0 \in {\bf Z}$. Moreover note that, by Assumption (v), there exist positive
 numbers $M_z$, $a_z$ and $d$ such that, for any ${\bf z}_0$ satisfying ${\rm dist}({\bf z}_0,{\cal
 A}_0) \le d$, the following bound holds
 \beeq{\label{zexp}
 {\rm dist}({\bf z}(t),{\cal
 A}_0) \leq M_z e^{-a_z t}\, {\rm dist}({\bf z}_0,{\cal
 A}_0)\,,
 }
 for all $t \geq 0$. From this and from the definition of
 $\nu(\cdot)$ (and in particular from the fact that $\nu(\cdot)$ is differentiable and
 vanishes on ${\cal A}_0$)  it follows that there is a constant $\gamma >0$ such that
 \[
 |\nu({\bf z}(t,{\bf z}_0))| \leq \gamma {\rm dist}({\bf z}(t),{\cal
 A}_0) \leq  \gamma M_z e^{-a_z t}\, {\rm dist}({\bf z}_0,{\cal
 A}_0)\,, \]
 for all $t \geq 0$ and all ${\bf z}_0$ satisfying ${\rm dist}({\bf z}_0,{\cal
 A}_0) \le d$.
 Consider now the candidate Lyapunov function $W({\bf x},t) =
 { {\bf x}}\tr \bar P(t) {\bf x}$,
 whose time derivative along the solution of (\ref{xssubs}) yields
 \[
 \ba{rcl}
 \dot W({\bf x},t) &=& -{{\bf x}}\tr Q(t) {\bf x} + 2 {{\bf x} }\tr
 \bar P(t) {\bf h}({\bf z}(t,{\bf z}_0),{\bf x}) + 2 {{\bf x}}\tr \bar P(t) \nu({\bf z}(t,{\bf z}_0))\\[2mm]
  &\leq& -c_3 |{\bf x}|^2 + 2 c_2 |{\bf x}||{\bf h}({\bf z}(t,{\bf z}_0),{\bf x})| +
  2 c_2 |{\bf x}| |\nu({\bf z}(t,{\bf z}_0))|\\[2mm]
  &\leq& -c_3 |{\bf x}|^2 + 2 c_2 |{\bf x}||{\bf h}({\bf z}(t,{\bf z}_0),{\bf x})| +
  2 c_2 |{\bf x}|\gamma M_z  e^{-a_z t}\, {\rm dist}({\bf z}_0,{\cal
 A}_0)\,.
  \ea
 \]
 Picking $\epsilon \leq \dst {c_3 \over 4 c_2}$ and
 $\delta_\epsilon$ accordingly, it follows that
 \[\ba{rcl}
 |{\bf x}(t)|\leq \delta_\epsilon \quad &\Rightarrow& \dst \quad \dot
 W({\bf x}(t),t) \leq -{c_3 \over 2} |{\bf x}(t)|^2 + 2 c_2 \gamma
 M_z \, |{\bf x}(t)| \, e^{-a_z t}\, {\rm dist}({\bf z}_0,{\cal
 A}_0)\\ \\
 \quad &\Rightarrow& \dst \quad \dot
 W({\bf x}(t),t) \leq -{c_3 \over 4 c_2}W({\bf x}(t),t) +
 4 {(c_2 \gamma M_z)^2 \over c_3} \,
 e^{-2a_z t}\, [{\rm dist}({\bf z}_0,{\cal
 A}_0)]^2\,.
 \ea\]
>From this, standard arguments can be invoked to claim the
existence of positive numbers $r_x < \delta_\epsilon$, $r_z<d$,
$A_x$, $A_y$, $\lambda_x$, ,$\lambda_z$ such that
 that if $|{\bf x}_0|<r_x$ and ${\rm dist}({\bf z}_0,{\cal
 A}_0) \leq
 r_z$ then the trajectory ${\bf x}(t)$ of (\ref{xssubs}) can be
 bounded as
 \[
 |{\bf x}(t)| \leq A_x e^{-\lambda_x t}|{\bf x}_0|+ A_z e^{-\lambda_z t}{\rm dist}({\bf z}_0,{\cal
 A}_0)\,.
 \]
 This proves the Lemma. \end{proof} $\triangleleft$

\section{Adaptive output regulation}

We return now to the closed loop system obtained from the
interconnection of (\ref{plard1}), (\ref{exo1}) and
(\ref{control}).
 As mentioned before, this system, viewed as a system with input
$v$ and output $e=e_1$ has relative degree 1. To put it in
``normal form", we use, instead of (\ref{newvar}), the change
variables \beeq{\label{newvar2}\ba{rcl}\tilde \theta &=& \hat
\theta - \theta(\rho) - \beta x
\\[2mm]
\eta &=& \xi - M[\hat \theta - \theta(\rho)] - Kx\,.\ea}

This, after some simple algebra and some obvious rearrangement of
terms, yields a system of the form
\beeq{\label{closedloop2}\ba{rcl}
 \dot \varrho &=& 0 \\[2mm]
 \dot w &=& s(\varrho,w) \\[2mm]
 \dot z &=&  f_0(\varrho,w,z)+f_1(\varrho, w,z,e)e \\[2mm]
 \dot \eta &=& A\eta +
 b\beta\tr \tilde \theta - K[q(\varrho,w,z,0)+\eta_1] +
 \phi(\eta_1) + \Omega(\eta_1)\theta +
 \delta_1(\varrho,w,z,e,X,\eta_1)\,e\\[2mm]
\dot {\tilde \theta} &=& -\beta [q(\varrho,w,z,0)+\eta_1] +
\delta_2(\varrho,w,z,e,X,\eta_1)\,e - \mbox{dzv}_\ell(\tilde
\theta + \theta(\varrho))\\[2mm]
\dot X &=& FX + G\Omega(\eta_1)+\delta_3(\eta_1,e)\,e\\
[2mm]\dot e &=&
-[q(\varrho,w,z,0)+\eta_1]+\vartheta(\varrho,w,z,e)e+v\,.\ea} in
which $\delta_1(\cdot)$, $\delta_2(\cdot)$, $\delta_3(\cdot)$ and
$\vartheta(\cdot)$ are continuously differentiable functions of
their arguments.

A more succinct form can be obtained setting $\bf w$ as in section
(\ref{SecPA}) and
 \[\ba{rcl} {\bf x} &=& {\rm col}(\eta,\tilde \theta, X_1,\ldots X_q)\,, \ea\]
(where $X_i$ denotes the $i$-th column of $X$) in which case, the
system in question can be rewritten in the form\,\footnote{With a
minor abuse of notation we have replaced $f_0(\varrho,w,z)$ and
$f_1(\varrho,w,z,e)$ by $f_0({\bf w},z)$ and, respectively,
$f_1({\bf w}, z,e)$.}
 \beeq{\label{closbold}
 \ba{rcl}
 \dot{\bf w} &=& {\bf s}({\bf w})\\[2mm]
 \dot z &=& f_0({\bf w},z) + f_1({\bf w},z,e)e\\[2mm]
 \dot {\bf x} &=& {\bf g}_0({\bf w},z,{\bf x})+ {\bf g}_1({\bf w},z,{\bf
 x},e)e
 \\[2mm]\dot e &=& {\bf h}({\bf w},z,{\bf x})+ {\bf k}({\bf w},z,{\bf x}, e)e+v\,.
 \ea}
 In this notation, the set of equations
 \beeq{\label{zerodynbold}\ba{rcl}
 \dot{\bf w} &=& {\bf s}({\bf w})\\[2mm]
 \dot z &=& f_0({\bf w},z)\\[2mm]
 \dot {\bf x} &=& {\bf g}_0({\bf w},z,{\bf x})
 \ea}
 is a succinct version for the set of equations (\ref{zerodyn}),
whose asymptotic properties have been analyzed in the previous
section. More precisely, under the hypotheses introduced earlier,
the positive orbit of ${\bold  Z}\times  \Xi \times \Theta \times
{\bf X}$ under the flow (\ref{zerodynbold}) is bounded and all
trajectories uniformly converge to the compact invariant set
$\omega({\bf B})$ described in Proposition \ref{PR1}. Moreover,
the function ${\bf h}({\bf w},z,{\bf x})$, which is a succinct
version of the quantity $-[q(\varrho,w,z,0 )+\eta_1]$ in
(\ref{closedloop2}), vanishes on the set $\omega({\bf B})$. With
this in mind we are now in the position to formulate the final
result of the paper which states that the controller
(\ref{control}) completed with
 \beeq{\label{vgx}
 v = - k e
 }
 solves the problem of output regulation if $k$ is chosen sufficiently large.

\begin{proposition} Consider system (\ref{plard1}) with exosystem (\ref{exo1}).
Let  ${\bf W},Z, E$ be fixed compact sets of initial conditions,
for which the assumptions (i)-(iv) indicated in section
\ref{SecPA} are supposed to hold. Suppose, in addition, that
assumption (v) introduced in section \ref{sec3} holds. Consider
the controller (\ref{control}) completed with (\ref{vgx}) and
initial conditions in a fixed compact set ${\bf K}$. Then, there
exists a number $k^\ast>0$ such that if $k\ge k^\ast$ the positive
orbit of ${\bf W} \times Z \times  E \times {\bf K}$ in the closed
loop system is bounded and $e(t)\to 0$ as $t\to \infty$.
\end{proposition}

\begin{proof}
 The result directly follows from Proposition \ref{PropApp} of
 Appendix \ref{AppA}. In particular it is easy to check that system
 (\ref{closbold}) -- (\ref{vgx}) can be viewed as a system of the form (\ref{sysappA}),
 the role of $x$ in (\ref{sysappA}) being played here by the
 one-dimensional variable $e$. The properties established for
 (\ref{zerodynbold}) and the fact that {\bf h}({\bf w},z,{\bf x}) vanishes on $\omega({\bf B})$
 show that all the assumptions of Proposition \ref{PropApp}
 are satisfied. Thus, the desired result follows by taking a large value of $k$.\end{proof}$\triangleleft$

 \begin{remark} The previous Proposition indicates that the proposed controller
 (\ref{control}) completed with (\ref{vgx}) solves the problem of output regulation
 for the relative degree one system (\ref{plard1}). Bearing in
 mind the
 discussion at the end of section \ref{SecPA}, though, it follows that a controller
 of the form (\ref{control}), completed with
 \[ v = - k(e_r + g^{r-1} a_0 e_1 + g^{r-2} a_1 e_2 + \ldots + g
 a_{r-2} e_{r-1})\,,\]
is able to solve the problem of output regulation for the original
plant (\ref{pla1}), if $g$ is large enough. In this respect, it is
worth stressing that the assumptions under which the proposed
controller solves the problem need only to be checked on the
original system (\ref{pla1}) and not necessarily on the
transformed, relative degree one, system (\ref{plaaux}). As a
matter of fact, we have already shown, in Lemma \ref{lemmard},
that assumptions (i) through (iv) on system (\ref{pla1}) imply
identical properties on system (\ref{plaaux}). For the sake of
coherence, it remains to show that the fulfilment of assumption
(v) on system (\ref{pla1}) implies the fulfillment of the
corresponding assumption on system (\ref{plaaux}). But this is a
trivial matter, in view of the fact that the assumption in
question is determined (once the matrices $A$, $C$, $F$, $G$ and
the map $\Omega(\cdot)$ have been fixed) only by the restriction
of $\tau_1({\bf z})$ to the invariant set ${\cal A}_0$. As shown
in the proof Lemma \ref{lemmard}, the map $\tilde \tau (\tilde
{\bf z})$ which makes assumption (iv) satisfied for (\ref{plaaux})
is such that $\tilde \tau (\tilde {\bf z})=\tau({\bf z})$ for
$\tilde {\bf z}=({\bf z},0)$ and $ {\bf z} \in {\cal A}_0$, and
therefore, if system (\ref{pla1}) has the property (v), an
identical property holds for the transformed system
(\ref{plaaux}). $\triangleleft$

 \end{remark}

\bigskip
\begin{center} {\bf \Large Appendix}
\end{center}
\appendix

\bigskip

 \section{A small-gain property}\label{AppA}

 Consider a system of the form
 \beeq{\label{sysappA}
 \ba{rcl}
 \dot \varrho &=& 0\\
 \dot {w} &=& {s}(\varrho, {w})\\
 \dot z &=& f_0(\varrho, w,z) + \ell(\varrho, w,z,x)\\
 \dot x &=& q_0(\varrho,  w,z) + r(\varrho, {w},z,x)  + g A x
 \ea
 }
 in which $(\varrho, w, z, x) \in \Real^p \times \Real^s \times \Real^n \times
 \Real^m$. Let the  functions $f_0(\cdot), q_0(\cdot), g(\cdot), \ell(\cdot), s(\cdot)$
be continuously differentiable and, moreover, let
 $\ell(\varrho, {w},z,0)=0$ and $r(\varrho, {w},z,0)=0$ for
 all $(\varrho, {w},z) \in \Real^p \times \Real^s \times \Real^n$.
 $A$ is a given Hurwitz matrix and $g$ is a positive
 number. As in section
 \ref{SecPA}, let ${\cal P} \subset \Real^p$, $W \subset \Real^s$,
 $Z \subset \Real^n$ denote compact sets of
 initial conditions for $\varrho$, $w$, $z$,
 set ${\bf z} = \mbox{col}(\varrho, w, z)$
 and ${\bf Z} = {\cal P} \times W \times Z$.
 Suppose that the autonomous system
 \beeq{ \label{zdapp}
 \ba{rcl}
 \dot \varrho &=& 0\\
 \dot {w} &=& {s}(\varrho, {w})\\
 \dot z &=& f_0(\varrho, w,z)\,,
 \ea
 }
 with initial conditions in the compact set ${\bf Z}$,
 satisfies assumptions (i), (ii), (iii) of section
 \ref{SecPA} and, coherently with the abbreviated notation used throughout the paper, set  ${\cal A}_0= \omega({\bf Z})$.
 The following lemma presents describes some relevant properties of (\ref{sysappA}),
  proven
 \cite{BIP}, which instrumental in proving the desired results.

 \begin{lemma} \label{lemmaAppB}
 Consider system (\ref{sysappA}) under the assumptions specified
 above,
with initial conditions in ${\bf Z} \times X$ with
 $X\subset \Real^m$ a compact set. Then the following holds:

 \medskip

 \noindent (a) there exist positive numbers $d$, $M$, $a$ and $\gamma$
 such that if
 \[
 \mbox{\em dist}({\bf z}(0), {\cal A}_0) \leq d \qquad \mbox{and}
 \qquad |x(t)| \leq d \quad \mbox{for all }  t \geq 0
 \]
  then
 \[
 \mbox{\em dist}({\bf z}(t), {\cal A}_0)
 \leq M e^{- a t} \mbox{\em dist}({\bf z}(0), {\cal A}_0)
 + \gamma \dst \max_{\tau \in [0,t]} |x(\tau)|
  \qquad \mbox{for all }  t \geq 0.\]

 \medskip

 \noindent (b) for all $\epsilon>0$ there exist $g^\star>0$ and $T>0$ such
 that for all $g \geq g^\star$ the positive orbit of ${\bf Z} \times X$ under
 the flow of (\ref{sysappA}) is bounded and
 \[
 \mbox{\em dist}({\bf z}(t), {\cal A}_0) \leq
 \epsilon\,, \qquad
 |x(t)| \leq \epsilon \qquad \mbox{for all }  t \geq T.
 \]
 \end{lemma}

 The previous lemma provides the tools needed to study  the
 asymptotic behavior of the system (\ref{sysappA})
 under the additional hypothesis that the function $q_0(\varrho, w, z)$
 vanishes on ${\cal A}_0$ (or, what is the same, that ${\cal A}_0 \times \{0\}$ is invariant for
 (\ref{sysappA})).  This is specified in
 the next proposition.

 \begin{proposition}\label{PropApp}
 Consider system (\ref{sysappA}) under the assumptions specified
 above and assume, in addition,  that $q_0(\varrho,w,z)= 0$
 for all $(\varrho,w,z) \in {\cal A}_0$. Then for any compact
 set $X$ there exists $g_1^\star>0$ such that, for
 all $g \geq g_1^\star$, the positive orbit of
 ${\bf Z} \times {X}$ under the flow of (\ref{sysappA}) is
 bounded and $\dst \lim_{t \rightarrow \infty} x(t) =0$.
 \end{proposition}

 \begin{proof}
 The proof is an easy consequence of the results of Lemma
 \ref{lemmaAppB} and of the small gain theorem.
 As a matter of fact, pick $\epsilon\leq d$ and set
 \[
 {\cal S}_\epsilon = \{ ({\bf z},x)\in \Real^{p+s+n} \times \Real^m \; : \;
 \mbox{dist}({\bf z},{\cal A}_0) \leq \epsilon\,, \quad |x|
 \leq \epsilon \}
 \]
 From property (b), it is seen that if $g\geq g^\star$,
 any initial condition in ${\bf Z} \times
 X$ produces a trajectory  of (\ref{sysappA}) which is bounded in
 forward time and satisfies $({\bf z}(t),x(t))\in {\cal S}_\epsilon$
 for all $t \geq
 T$.
 From property (a), it is seen that
 \[
 \mbox{dist}({\bf z}(t-T),{\cal A}_0) \leq M e^{-a(t-T)}
 \mbox{dist}({\bf z}(T),{\cal A}_0) + \gamma \dst
 \max_{\tau \in [T,t-T]} |x(\tau)|
 \]
 for all $t \geq
 T$.
 Note that the differentiable function $r(\varrho,w,z,x)$, which vanishes for $x=0$,
 can be estimated as
 \[
 |r(\varrho,w,z,x)| \le   \alpha |x|
 \]
 for all $({\bf z},x) \in {\cal S}_\varepsilon$, while the differentiable function $q_0(\varrho,w,z)$ which vanishes on ${\cal A}_0$, can be estimated as
 \[
 |q_0(\varrho,w,z)| \leq \beta \, \mbox{dist}({\bf z},{\cal A}_0)
 \]
 for all ${\bf z} \in \cal S_\varepsilon$. Now, let $P>0$ denote the
 solution of $PA + A\tr P = -I$ and by $\underline \lambda$ and $\bar \lambda$
 respectively the smallest and largest eigenvalue of $P$.
 Standard arguments can be used to show that, for all $t \geq T$,
 \[
 |x(t)| \leq \sqrt{\bar \lambda \over \underline \lambda} \,
 e^{- {\lambda_g \over 4 \bar \lambda} (t-T)} |x(T)|
  + {4 \beta \over \lambda _g} |P| \, \max_{\tau \in
  [T,t-T]} \, \mbox{dist}({\bf z},{\cal A}_0)
 \]
 where $\lambda_g = g - 2 |P| \alpha$.
 Hence the result follows by classical small gain arguments if
 $g_1^\star \geq g^\star$ is picked so that  the small gain condition
 \[
 \lambda_{g_1^\star} > 4 \beta |P| \gamma
 \]
 is fulfilled. This completes the proof of Proposition \ref{PropApp}.
 $\triangleleft$
\end{proof}

\bigskip
\begin{center} {\bf \Large Acknowledgements}
\end{center}
\medskip
The authors wish to thank Christopher I. Byrnes and Laurent Praly
for fruitful discussion and helpful suggestions during the
preparation of the paper.


\begin{thebibliography}{99}


\bibitem{BaGe} G. Bastin and M. R. Gevers, Stable adaptive
observers for  non-linear time varying systems, {\em IEEE Trans.
Autom. Contr.\/}, {\bf AC-33}: 650--657, 1988.

\bibitem{BI03} C.I. Byrnes and A. Isidori, Limit sets, zero dynamics and internal models
in the problem of nonlinear output regulation, {\em IEEE Trans. on
Automatic Control}, {\bf AC-48}, pp. 1712--1723, 2003.

\bibitem{BI03bis} C.I. Byrnes and A. Isidori, Nonlinear Internal Models for Output
Regulation, Preprint arXiv: math.OC/0311223.

\bibitem{BIP} C.I. Byrnes, A. Isidori and L. Praly, On the Asymptotic Properties of a System Arising in
Non-equilibrium Theory of Output Regulation, Preprint of the
Mittag-Leffler Institute, Stockholm, 18, 2002-2003, spring.

\bibitem{BDIK}
C.I. Byrnes, F. Delli Priscoli, A. Isidori and W. Kang,
\newblock Structurally stable output regulation of nonlinear systems.
\newblock {\em Automatica}, {\bf 33}: 369--385, 1997.

\bibitem{Da} E.J. Davison, \newblock The robust control of a servomechanism problem for
linear time-invariant multivariable systems, \newblock {\em IEEE
Trans. Autom. Contr.}, {\bf AC-21}: 25--34, 1976.

\bibitem{Fr} B.A. Francis, \newblock The linear multivariable regulator problem, \newblock
{\em SIAM J. Contr. Optimiz.}, {\bf 14}: 486--505, 1977.


\bibitem{FW76}
B.A. Francis and W.~M. Wonham.
\newblock The internal model principle of control theory.
\newblock {\em Automatica}, {\bf 12}: 457--465, 1976.

\bibitem{Hale91} J.K. Hale, L.T. Magalh\~aes and W.M. Oliva, {\em
Dynamics in Infinite Dimensions}, Springer Verlag (New York, NY),
2002.

\bibitem{HL}
J. Huang and C.F. Lin.
\newblock On a robust nonlinear multivariable servomechanism problem.
\newblock {\em IEEE Trans. Autom. Contr.\/}, {\bf AC-39}: 1510--1513, 1994.

\bibitem{HR}
J. Huang and W.J. Rugh.
\newblock On a nonlinear multivariable servomechanism problem.
\newblock {\em Automatica}, {\bf 26}:963--972, 1990.

\bibitem{IA} A. Isidori, \newblock A remark on the problem of semiglobal
nonlinear output regulation, {\em IEEE Trans. on Automatic
Control}, {\bf AC-42}: 1734-1738, 1997.

\bibitem{IB}
A. Isidori and C.I. Byrnes.
\newblock Output regulation of nonlinear systems.
\newblock {\em IEEE Trans. Autom. Contr.\/}, {\bf AC-25}: 131--140, 1990.

\bibitem{Kh2} H. Khalil, Robust servomechanism output feedback controllers for
feedback linearizable systems, \newblock {\em Automatica}, {\bf
30}: 1587--1599, 1994.

\bibitem{Khalil} H.K. Khalil, {\em Nonlinear Systems},
                Prentice Hall,2nd edition, Upper Saddle River, NJ,
                1996.

\bibitem{MaToBook} R. Marino and P. Tomei, {\em Nonlinear Control Design: Geometric, Adaptive,
\& Robust}, Prentice Hall (New York), 1995.

\bibitem{MaTo92} R. Marino and P. Tomei, Global adaptive observers
for nonlinear systems via filtered transformations, {\em IEEE
Trans. on Automatic Control}, {\bf AC-37}, pp. 1239--1245, 1992.

\bibitem{Rugh} W.J. Rugh, {\em Linear System Theory},
Prentice Hall (New York), 1996.

\bibitem{SIM99a} A. Serrani, A. Isidori and L. Marconi,
 Semiglobal output regulation for minimum-phase systems,
 {\em International Journal on Robust and Nonlinear Control}, {\bf 10}, pp.
379--396, 2000.

\bibitem{SIM} A. Serrani, A. Isidori and L. Marconi,
Semiglobal nonlinear output regulation with adaptive internal
model, {\em IEEE Trans. Autom. Contr.}, {\bf AC-46}: 1178-1194,
2001.

\bibitem{JS} J. Szarski, {\em Differential Inequalities}, Polska
Akademia Nauk (Warszawa), 1967.

\bibitem{TP}
A.R. Teel and L. Praly,
\newblock Tools for semiglobal stabilization by partial
state and output feedback.
\newblock {\em SIAM J. Control Optim.}, {\bf 33}, pp. 1443--1485, 1995.

\bibitem{Wi} F.W. Wilson, Smoothing derivatives of functions and
applications, {\em Trans. Amer. Math. Soc.}, {\bf 139}: 413--428,
.

\bibitem{Yoshizawa} T. Yoshizawa, {\em Stability Theory and the
Existence of Periodic Solutions and Almost Periodic Solutions},
Springer Verlag (New York, NY), 1975.

\end{thebibliography}
\end{document}